\documentclass[final]{siamltex}
\pdfoutput=1

\usepackage{graphicx}

\usepackage{amsmath} 
\usepackage{amssymb}  
\usepackage{color}
\usepackage{graphics}
\usepackage{epsfig}
\usepackage[ruled,vlined]{algorithm2e}
\usepackage{comment}
\usepackage{dsfont}

\newtheorem{remark}{Remark}[section]
\newtheorem{assumption}{Assumption}[section]
\newcommand{\beq}{\begin{equation}}
\newcommand{\eeq}{\end{equation}}

\newcommand{\C}{\mathbb C}

\newcommand{\R}{\mathbb R}

\def\iu{\imagunit}
\def\e{{\rm e}}

\def\eps{\varepsilon}

\def\l{\lambda}

\newcommand{\tr}{{\rm trace}}

\newcommand{\la}{\left\langle}
\newcommand{\ra}{\right\rangle}
\newcommand{\matlab}{{\sc Matlab }}

\newcommand{\Lameps}{\Lambda_\varepsilon}

\renewcommand{\Re}{{\mbox{\rm Re}}}
\renewcommand{\Im}{{\mbox{\rm Im}}}

\newcommand{\imagunit}{{\bf i}}

\newcommand{\goes}{\rightarrow}

\newcommand{\wt}{}
\newcommand{\Deltaopt}{\Delta_{\rm opt}}

\newcommand{\mussv}{{\tt mussv}}

\newcommand{\esm}{\end{smallmatrix}\right]}

\newcommand{\BP}{\mathbb B}
\newcommand{\BPR}{{\mathbb B}}
\newcommand{\BPC}{{{\mathbb B}^{*}}}
\newcommand{\BPCO}{{{\mathbb B}_1^*}}
\newcommand{\Id}{{\rm I}}
\newcommand{\One}{\mathds{1}}

\newcommand{\Z}{Z}
\newcommand{\F}{\mathbf{F}}
\newcommand{\cM}{\mathcal{M}}
\newcommand{\lbn}{\mu^{\ell}_{\rm \scriptscriptstyle New}}
\newcommand{\lbo}{\mu^{\ell}_{\rm \scriptscriptstyle PD}}
\newcommand{\z}{\zeta}
\newcommand{\lz}{\zeta}
 
\title{A novel iterative method to approximate structured singular values}

\author{Nicola Guglielmi\footnotemark[1] \and Mutti-Ur Rehman\footnotemark[2] \and Daniel Kressner\footnotemark[3]}

\begin{document}

\renewcommand{\thefootnote}{\fnsymbol{footnote}}
\footnotetext[1]{Dipartimento di Ingegneria Scienze Informatiche e Matematica (DISIM),
Universit\`a degli Studi di L' Aquila,
Via Vetoio - Loc.~Coppito,
and Gran Sasso Science Institute (GSSI), via Crispi 7,
I-67010    L' Aquila,  Italy. Email: {\tt guglielm@univaq.it}}

\footnotetext[2]{Gran Sasso Science Institute (GSSI), via Crispi 7,
I-67010    L' Aquila,  Italy. Email: {\tt mutti.abbasi@gssi.infn.it}}

\footnotetext[3]{EPFL-SB-MATHICSE-ANCHP,
Station 8,
CH-1015 Lausanne, Switzerland. Email: {\tt daniel.kressner@epfl.ch}}

\date{21 March 2016}

\maketitle

\begin{abstract}
A novel method for approximating structured singular values (also known as $\mu$-values) 
is proposed and investigated. These quantities constitute an important 
tool in the stability analysis of uncertain linear control systems as well as in 
structured eigenvalue perturbation theory. Our approach consists of an inner-outer iteration. In the outer iteration,
a Newton method is used to adjust the perturbation level. The inner iteration solves a gradient system associated with
an optimization problem on the manifold induced by the structure.
Numerical results and comparison with the well-known {\rm Matlab} 
function \mussv, implemented in the \matlab Control Toolbox, illustrate the behavior 
of the method. 
\end{abstract}



\section{Introduction}

The structured singular value (SSV)~\cite{PD93} is an important and versatile tool in control, as it allows to
address a central problem in the analysis and synthesis of control systems: To quantify the
stability of a closed-loop linear time-invariant systems subject to structured perturbations.
The class of structures addressed by the SSV is very general and allows to cover all types of parametric uncertainties that
can be incorporated into the control system via real or complex linear fractional transformations.
We refer to~\cite{BRQ98,CFN96a,CFN96b,HP10,Ka11,KHP06,
PD93,QBRDYD95,ZDG96} and the references therein for examples and applications of the SSV.

The versatility of the SSV comes at the expense of being notoriously hard, in fact NP hard~\cite{BYDM94}, to compute.
Algorithms used in practice thus aim at providing upper and lower bounds, often resulting in a coarse estimate of 
the exact value. An upper bound of the SSV provides sufficient conditions to 
guarantee robust stability, while a lower bound provides sufficient conditions 
for instability and often also allows to determine structured perturbations that 
destabilize the closed loop linear system.

The widely used function \mussv{} in the \matlab Control Toolbox 
computes an upper bound of the SSV using diagonal balancing / LMI techniques~\cite{YND92,FTD91}.
The lower bound is computed by a generalization of the power method developed in~\cite{YNDP94,PFD88}.
This algorithm {\em resembles a mixture of the power methods for computing the spectral radius and the largest singular value}, which is not
surprising, since the SSV can be viewed as a generalization
of both. 
When the algorithm converges, a lower bound of the SSV
results and this is always an equilibrium point of the iteration.
However, in contrast to the standard power method, there are, in general, several stable equilibrium points and not all of them correspond to the SSV.
In turn, one cannot guarantee convergence to the exact value but only to a lower bound. 
We remark that, despite this drawback, \mussv{} is a very reliable and powerful routine, which reflects
the state of the art in the approximation of the SSV.

In this paper, we present a new approach to computing a lower bound of the SSV 
associated with general mixed real/complex perturbations. The main ingredient of our new algorithm
is a gradient system that evolves perturbations on a certain
matrix manifold towards critical perturbations. Among the theoretical properties established for this gradient system, we prove a monotonicity property that indicates robustness and can also be exploited in the numerical discretization.
We show several numerical examples for which our algorithm provides tighter bounds than those computed by \mussv{}.

\subsection{Overview of the article}

Section \ref{sec:fw} provides the basic framework for the proposed methodology.
In particular, we explain how the computation of the SSV can be addressed by an inner-outer algorithm, where the
outer algorithm determines the perturbation level $\eps$ and the inner algorithm determines a (local) extremizer of the
structured spectral value set. Moreover, an example illustrates that the output \mussv{} may fail to satisfy a necessary
condition for optimality. 


In Section \ref{sec:compl} we develop the inner algorithm for the case of complex structured perturbations.
An important characterization of extremizers shows that we can restrict ourselves to a manifold of 
structured perturbations with normalized and low-rank blocks.
A gradient system for finding extremizers on this manifold is established and analyzed.

Section~\ref{sec:mix} extends the results of Section~\ref{sec:compl} to perturbations
with complex full blocks alternated and mixed complex/real repeated scalar blocks. 

The outer algorithm is addressed in Section~\ref{sec:fast}, where a Newton method for determining the correct perturbation level $\eps$ is developed.
The algorithm proposed in this work is presented in Section~\ref{sec:algoR}.

Finally, in Section \ref{sec:compres}, we present a range of numerical 
experiments to compare the quality of the lower bounds obtained with our algorithm to those obtained with \mussv{}.

%

\section{Framework}
\label{sec:fw}

We consider a matrix $M \in \C^{n\times n}$ and an underlying perturbation set
with prescribed block diagonal structure,
\begin{equation}
\label{eq:BP}
\BP = \big\{ \diag \left( \delta_1 \Id_{r_1}, \ldots, \delta_s \Id_{r_S},
\Delta_1, \ldots,\Delta_F \right), \delta_i \in \C \mbox{($\R$)},
\Delta_j \in \C^{m_j \times m_j} \mbox{($\R^{m_j \times m_j}$)} \big\},
\end{equation}
where $\Id_{r_i}$ denotes the $r_i \times r_i$ identity matrix.
Each of the scalars $\delta_i$ and the $m_j \times m_j$ matrices $\Delta_j$
may be constrained to stay real in the definition of $\BP$.
The integer $S$ denotes the number of repeated {\em scalar} blocks (that is, scalar multiples 
of the identity) and $F$ denotes the number of {\em full} blocks.
This implies $\sum_{i=1}^{S} r_i + \sum_{j=1}^{F} m_j = n$.
In order to distinguish complex and real scalar blocks, we assume
that the first $S'\le S$ blocks are complex while the (possibly) remaining $S-S'$ blocks 
are real. Similarly we assume that the first $F'\le F$ full blocks are complex 
and the (possibly) remaining $F-F'$ blocks are real. The literature (see, e.g., \cite{PD93})
usually does not consider real full blocks, that is, $F'=F$. In fact, in
control theory, full blocks arise from uncertainties associated to the frequency 
response of a system, which is complex-valued.  

For simplicity, we assume that all full blocks are square, although this
is not necessary and our method extends to the non-square case in a
straightforward way. Similarly, the chosen ordering of blocks should not be viewed as a limiting assumption; it merely simplifies notation.  

The following definition is given in \cite{PD93}, where $\|\cdot\|_2$ denotes
the matrix $2$-norm and $\Id$ the $n\times n$ identity matrix.  
\begin{definition} \label{def:mu}
Let $M \in \C^{n\times n}$ and consider a set $\BP$ of the form~\eqref{eq:BP}. Then the SSV (or $\mu$-value) $\mu_\BP(M)$ is defined as
\begin{equation}
\mu_\BP(M)  := \frac{1}{\min\left\{ \|\Delta\|_2 : \Delta \in \BP, 
\det(\Id - M \Delta ) = 0 \right\}}.
\label{eq:muB}
\end{equation}
\end{definition}
In Definition~\eqref{def:mu} and in the following, we use the convention that the minimum over an empty set is $+\infty$. In particular,
$\mu_\BP(M) = 0$ if $\det(\Id - M \Delta ) \neq 0$ for all $\Delta \in \BP$.

Note that $\mu_\Delta$ is a positively homogeneous function, i.e.,
\[
\mu_\BP(\alpha M) = \alpha \mu_\BP(M) \qquad \mbox{for any $\alpha \ge 0$}.
\]
For $\BP = \C^{n\times n}$, it follows directly from Definition~\ref{def:mu} that $\mu_\BP(M) = \|M\|_2$. For general $\BP$, the SSV can only become smaller and we thus have the upper bound $\mu_\BP(M) \le \|M\|_2$. This can be refined further by exploiting the properties of $\mu_\BP$, see \cite{ZDG96}.
These relations between $\mu_\BP(M)$ and $\|M\|_2$, the largest singular value of $M$, justifies the name \emph{structured singular value} for $\mu_\BP(M)$.

The important special case when $\BP$ only allows for complex perturbations, that is, $S=S'$ and $F=F'$,
deserves particular attention. 
In this case we will write $\BPC$ instead of $\BP$. 
Note that $\Delta \in \BPC$ implies $\e^{\iu \varphi} \Delta \in \BPC$ for any $\varphi \in \R$. 
In turn, there is $\Delta \in \BPC$ such that $\rho(M \Delta)=1$ if and only if there is $\Delta' \in \BPC$, with
the same norm, such that $M \Delta'$ has the eigenvalue $1$, which implies $\det(\Id - M \Delta')=0$.
This gives the following alternative expression:
\begin{equation}
\mu_\BPC(M)  =  \frac{1}{\min\left\{ \|\Delta\|_2 : \Delta \in \BPC, 
\rho( M \Delta ) = 1 \right\}},
\label{eq:muBC}
\end{equation}
where $\rho(\cdot)$ denotes the spectral radius of a matrix.
For any nonzero eigenvalue $\lambda$ of $M$, the matrix $\Delta = \lambda^{-1} I$ satisfies the constraints of the minimization problem in~\eqref{eq:muBC}.
This establishes the lower bound  $\rho(M) \le \mu_\BPC(M)$ for the case of purely complex perturbations.
Note that $\mu_\BPC(M) = \rho(M)$ for $\BPC = \{ \delta \Id:\, \delta \in \C\}$.
Hence, both the spectral radius and the matrix 2-norm are included as (trivial) special cases of the SSV.

\subsection{A motivating example}
\label{sec:motex}


Consider the $3 \times 3$ matrix 
\begin{eqnarray*}
M & = & \left( \begin{array}{rrr}
-1 + \iu &   1 - \iu  & -1 + \iu \\
-1 + \iu &  -1  & \iu \\
\iu &  -1 - \iu & 1 - \iu	
\end{array} \right),
\end{eqnarray*}
where $\iu$ denotes the imaginary unit, along with the perturbation set
\[
\BP = \big\{ \diag( \delta_1 \Id_{2}, \Delta_1 ):\, \delta_1 \in \R, \ \Delta_1 \in \C^{1,1} \big\}.
\]
Applying the {\sc Matlab} function \mussv\footnote{In all experiments we have used \mussv{} with its default
parameters.}
yields the bounds
\begin{equation}
0.9807\ldots \le \mu_\BP(M) \le 2.2477\ldots.
\label{eq:estexil2}
\end{equation}
The large difference between the lower and upper bounds is caused by the lower bound.
The perturbation determining the lower bound is given by  $\widehat\eps \widehat\Delta$  with
\begin{eqnarray*}
\widehat\Delta =
\left( \begin{array}{rrr}
 -0.368473881\ldots & 0 & 0 \\                     
  0 &           -0.368473881\ldots    & 0 \\                      
  0 &            0               & -0.673755352\ldots - 0.738954481\ldots\iu
\end{array} \right)
\end{eqnarray*}
and $\widehat\eps = 1.019727084\ldots$.
The scaling has been chosen such that
$\| \widehat\Delta \|_2 = 1$. However, not all blocks of $\widehat\Delta$ have unit norm; the 
$2 \times 2$ repeated scalar block of $\widehat\Delta$ has norm $0.368473881\ldots$. 
We will see in Theorem \ref{th:maxmix} below that this violates a necessary optimality condition for an extremizer $\Delta \in \BP$, 
which states that the spectral norm of all blocks of a normalized extremizer, under suitable conditions which are fulfilled here, should be one.

	
Applying our new algorithm, Algorithm~\ref{algo-ode} below, we obtain the 
perturbation $\eps^\star \Delta^\star$ with  
\begin{eqnarray*}
\Delta^\star =
\left( \begin{array}{rrr}
 -1 &  0 & 0 \\                     
  0 & -1 & 0 \\                      
  0 &  0 &  -0.989237164 - 0.146320991\ldots \iu 
\end{array} \right).
\end{eqnarray*}
and $\eps^\star = 0.445238645\ldots$, determining the lower bound \[\mu_\BP(M) \ge \lbn = 2.2459865301\ldots,\] which makes the estimate \eqref{eq:estexil2} substantially sharper. Note that both blocks of $\Delta^\star$ have unit norm.


\subsection{A reformulation based on structured spectral value sets}
\label{sec:reform}
The structured 
spectral value set of $M \in \C^{n\times n}$ with respect to a perturbation level $\eps$ is defined as
\begin{equation}
\Lambda_\eps^{\BP}(M) = 
\{\, \lambda \in \Lambda( \eps M \Delta ): \,
\Delta \in \BP, \ \| \Delta \|_2 \le 1 \}, 
\label{eq:svs}
\end{equation}
where $\Lambda(\cdot )$ denotes the spectrum of a matrix. 
Note that for purely complex $\BP^*$, the set~\eqref{eq:svs} is simply a disk centered at $0$. 
The set 
\begin{equation}
\Sigma_\eps^{\BP}(M) = 
\{\, \z = 1-\lambda:
\ \lambda \in \Lambda_\eps^{\BP}(M) \} 
\label{eq:ssvs}
\end{equation}
allows us to express the SSV defined in~\eqref{eq:muB} as 
\[
\mu_\BP(M) = \frac{1}{\arg \min\limits_{\eps>0} \left\{ 0 \in \Sigma_\eps^{\BP}(M) \right\}},
\]
that is, as a structured distance to singularity problem.
We have that $0 \not\in \Sigma_\eps^{\BP}(M)$ if and only if 
$\mu_\BP(M) < 1/\eps$. 

For a purely complex perturbation set $\BPC$, we can use~\eqref{eq:muBC} to alternatively express the SSV as 
\begin{equation} \label{eq:mucompl}
\mu_\BPC(M) = \frac{1}{\arg \min\limits_{\eps>0} 
\Big\{ \max\limits_{\lambda \in \Lambda_\eps^{\BPC}\!(M)} |\lambda| = 1 \Big\}}.
\end{equation}
We have that $\Lambda_\eps^{\BPC}(M) \subset D$, where $D$ denotes the open complex unit disk, if and only if 
$\mu_\BPC(M) < 1/\eps$.  

\subsection{Overview of the proposed methodology}
\label{sec:meth}

Let us consider the minimization problem
\begin{eqnarray}
\z(\eps) & = & \arg\min\limits_{\z \in \Sigma_\eps^{\BP}(M)} |\z|
\label{eq:cr}
\end{eqnarray}
for some fixed $\eps > 0$. By the discussion above, the SSV $\mu_\BP(M)$ is the reciprocal of the smallest value of $\eps$ for which $\z(\eps) = 0$.
This suggests a two-level algorithm: In the inner algorithm, we attempt to solve~\eqref{eq:cr}.
In the outer algorithm, we vary $\eps$ by an iterative procedure which exploits the knowledge of the exact derivative 
of an extremizer -- say $\Delta(\eps)$ -- with respect to $\eps$. We address~\eqref{eq:cr} by solving a system of ODEs.
In general, this only yields a local minimum of~\eqref{eq:cr} which, in turn, gives an upper bound for $\eps$ and hence a lower bound for $\mu_\BP(M)$.
Due to the lack of global optimality criteria for~\eqref{eq:cr}, the only way to increase the robustness of the method is to compute
several local optima.

The case of a purely complex perturbation set $\BPC$ can be addressed analogously by letting
the inner algorithm determine local optima for
\begin{equation}
\lambda(\eps)  =  \arg\max\limits_{\lambda \in \Lambda_\eps^{\BPC}(M)} |\lambda|,
\label{eq:crcomplex}
\end{equation}
which then yields a lower bound for $\mu_\BPC(M)$.

\section{Purely complex perturbations}
\label{sec:compl}

In this section, we consider the solution of the inner problem~\eqref{eq:crcomplex} in the estimation of $\mu_\BPC(M)$ for
$M \in \C^{n\times n}$ and a purely complex perturbation set
\begin{equation}
\label{eq:BPC} \nonumber
\BPC = \big\{ \diag( \delta_1 \Id_{r_1}, \ldots, \delta_S \Id_{r_S},
\Delta_1, \ldots,\Delta_F ):\, \delta_i \in \C,
\Delta_j \in \C^{m_j \times m_j} \big\}.
\end{equation}

\subsection{Extremizers}
\label{sec:extrc}

We will make use of the following standard eigenvalue perturbation result, see, e.g., \cite[Section II.1.1]{Kat95}. 
Here and in the following, we denote $\dot{\phantom{a}}= d/dt$.

\begin{lemma} \label{lem:eigderiv} 
Consider a smooth matrix family $C:\R \to \C^{n\times n}$ and let 
$\lambda(t)$ be an eigenvalue of $C(t)$
converging to a simple eigenvalue $\lambda_0$ of $C_0=C(0)$ as 
$t \goes 0$.  Then 
$\lambda(t)$ is analytic 
near $t=0$ with
$$
    \dot \lambda(0) = \frac{y_0^* C_1 x_0}{y_0^*x_0},
$$
where $C_1 = \dot{C}(0)$ and $x_0,y_0$ are right and 
left eigenvectors of $C_0$ associated to $\lambda_0$, that is,
$(C_0 -\lambda_0 \Id)x_0 = 0$ and $y_0^*(C_0-\lambda_0\Id)=0$.
\end{lemma}

Our goal is to solve the maximization problem~\eqref{eq:crcomplex}, which requires finding a perturbation $\Deltaopt$ such that
$\rho(\eps M \Deltaopt)$ is maximal among all $\Delta \in \BPC$ with $\| \Delta \|_2 \le 1$.
In the following, we call $\lambda$ a largest eigenvalue if $|\lambda|$ equals the spectral radius.
\begin{definition}
A matrix $\Delta \in \BPC$ such that $\| \Delta \|_2 \le 1$ and $\eps M \Delta$ has
a largest eigenvalue that locally maximizes the modulus of $\Lameps^{\BPC}(M)$
is called a {\em local extremizer}.
\end{definition}


The following result provides an important characterization of local extremizers.

\begin{theorem} \label{th:max}
Let 
$$
\wt\Deltaopt = \diag \left( \wt\delta_1 \Id_{r_1}, \ldots, \wt\delta_s \Id_{r_S},
\wt\Delta_1, \ldots,\wt\Delta_F \right), \quad \| \wt\Deltaopt \|_2 = 1, 
$$
be a local extremizer of $\Lambda_\eps^{\BPC}(M)$. We assume that $\eps M \wt\Deltaopt$
has a simple largest eigenvalue $\lambda  = |\lambda| \e^{\imagunit \theta}$,
with the right and left eigenvectors $x$ and $y$ scaled such that $s = \e^{\imagunit \theta} y^* x > 0$. 
Partitioning 
\begin{equation} \label{eq:partitionxz}
x  =  \left( x_1^{\rm T} \ \ldots \ x^{\rm T}_S, \ x^{\rm T}_{S+1} \ \ldots \ x_{S+F}^{\rm T} \right)^{\rm T}, \quad 
z  =  M^* y = \left( z_1^{\rm T} \ \ldots \ z^{\rm T}_S, \ z^{\rm T}_{S+1} \ \ldots \ z_{S+F}^{\rm T} \right)^{\rm T},
\end{equation} 
such that the size of the components $x_k, z_k$ equals the size of the $k$th block in $\wt\Deltaopt$, 
we additionally assume that
\begin{eqnarray}
&& z_k^* x_k \neq 0 \quad \forall \ k=1,\ldots,S 
\label{eq:A1}
\\[2mm]
&& \|z_{S+h}\|_2 \cdot \|x_{S+h}\|_2 \ne 0 \quad \forall \ h=1,\ldots,F. 
\label{eq:A2}
\end{eqnarray}
Then
$$
| \wt\delta_{k} | = 1 \quad \forall \ k=1,\ldots,S
\quad \mbox{and} \quad \| \wt{\Delta}_{h} \|_2 = 1 \quad \forall h=1,\ldots,F,
$$
that is, all blocks of $\wt\Deltaopt$ have unit $2$-norm.
\end{theorem}
\begin{proof}
\rm 
The result is proved by contradiction.
We first assume that $\| \wt{\Delta}_{h} \|_2 < 1$ for some $1 \le h \le F$ and consider the matrix-valued function
\begin{eqnarray}
\Delta(t) & = & \diag\big( \wt\delta_{1} \Id_{r_1},\ldots,\wt\delta_{S} \Id_{r_S},\wt{\Delta}_{1},
\ldots, \wt{\Delta}_{h} + t\, z_{S+h}\, x_{S+h}^*, \ldots, \wt{\Delta}_{F} \big),
\label{eq:Delta2}
\end{eqnarray}
which satisfies $\Delta(0)=\wt\Deltaopt$ and $\| \Delta(t) \|_2 \le 1$ 
for $t$ sufficiently small. 
Since $\lambda(0)=\lambda$ is simple, we can apply Lemma~\ref{lem:eigderiv} to $\eps \Delta(t) M$ and obtain
\begin{eqnarray}
\frac{d}{d t} |\lambda(t)|^2 \Big|_{t=0} & = &
2\,\Re(\overline\lambda \dot\lambda) =
2\,\Re\Bigl(\overline\lambda\, \eps \,\frac{y^* M \dot \Delta x} {y^*x} \Bigr)
\nonumber
\\
& = & 
2 \eps |\lambda| \Re\Bigl(\frac{y^* M \dot \Delta x} {\e^{\imagunit \theta} y^*x} \Bigr) =
2 \eps \frac{|\lambda|}{s} \Re(z^* \dot \Delta x  ).
\label{eq:lambda2}
\end{eqnarray} 
Inserting~\eqref{eq:Delta2} and exploiting \eqref{eq:A2}, we obtain
\begin{eqnarray*}
\frac{d}{d t} |\lambda(t)|^2 \Big|_{t=0} & = & 
2 \eps \frac{|\lambda|}{s} \| z_{S+h} \|_2^2 \cdot \| x_{S+h} \|_2^2 > 0,
\end{eqnarray*}
which contradicts the extremality of $|\lambda|$.
\smallskip

Let us now assume that $| \wt\delta_{k} | < 1$ for some $1 \le k \le S$ and
consider the matrix valued function
\begin{equation}
\Delta(t) = \diag \left( \wt\delta_{1} \Id_1,\ldots,\wt\delta_{k} \Id_k + 
t\, x_k^*z_k \Id_{k}, 
\ldots, \delta_{S} \Id_S, \wt\Delta_1, \ldots,\wt\Delta_F \right)
\nonumber \label{eq:Delta2b}
\end{equation}
which again satisfies $\Delta(0)=\wt\Deltaopt$ and $\| \Delta(t) \|_2 \le 1$ for $t$ sufficiently small. 
In analogy to the first part, Assumption~\eqref{eq:A1} implies
\begin{eqnarray*}
\frac{d}{d t} |\lambda(t)|^2 \Big|_{t=0} & = & 
2 \eps \frac{|\lambda|}{s} | z_k^* x_k |^2 > 0.
\end{eqnarray*}
This again gives a contradiction.
\end{proof}
\smallskip
\begin{remark}
Note that Assumptions \eqref{eq:A1} and \eqref{eq:A2} as well
as the simplicity of $\lambda$ are generic and commonly found in the literature on algorithms for the SSV, see, e.g.,~\cite[Sec. 7.2]{PD93}.
\end{remark}

The following theorem allows us to replace the full blocks in a local extremizer by rank-1 matrices.
\begin{theorem}
\label{th:max2}
Let 
$
\wt\Deltaopt = \diag \left( \wt\delta_1 \Id_{r_1}, \ldots, \wt\delta_s \Id_{r_S},
\wt\Delta_1, \ldots,\wt\Delta_F \right) 
$
be a local extremizer and let
$\lambda,x,z$ be defined and partitioned as in Theorem~\ref{th:max}.
Assuming that~\eqref{eq:A2} holds, every block $\wt\Delta_h$ has a singular value $1$ with associated 
singular vectors $u_h = \gamma_h z_{S+h}/\| z_{S+h}\|_2$  and $v_h = \gamma_h x_{S+h}/\| x_{S+h}\|_2$ 
for some $|\gamma_h| = 1$.
Moreover, the matrix 
$$
\Delta_* = \diag \left( \wt\delta_1 \Id_{r_1}, \ldots, \wt\delta_S \Id_{r_S},
u_1 v_1^*, \ldots,u_F v_F^* \right)
$$   
is also a local extremizer, i.e., $\rho(\eps M \Deltaopt) = \rho(\eps M \Delta_*)$.
\end{theorem}
\begin{proof}
Let $\hat z_{S+h} = z_{S+h}/\| z_{S+h}\|_2$, $\hat x_{S+h} = x_{S+h}/\| x_{S+h}\|_2$
and consider the matrix valued function
\begin{eqnarray}
\Delta(t) & = & \diag\left( \wt\delta_{1} \Id_1,\ldots,\wt\delta_{S} \Id_S,\wt{\Delta}_{1},
\ldots, (1-t) \wt{\Delta}_{h} + t\, \hat z_{S+h}\, \hat x_{S+h}^*, \ldots, \wt{\Delta}_{F} \right)
\nonumber \label{eq:Delta3}
\end{eqnarray}
which has $2$-norm bounded by $1$ for $t \in [0,1]$.
By Theorem~\ref{th:max},  $\| \wt\Delta_h \|_2 = 1$, which implies 
$|z_{S+h}^* \wt\Delta_h x_{S+h}| \le \| z_{S+h} \|\,\| x_{S+h} \|$.
Consequently,
\begin{eqnarray*}
\Re(z^* \dot \Delta x ) &=& \Re(z_{S+h}^* \wt\Delta_h x_{S+h} + z_{S+h}^* \hat z_{S+h}\, \hat x_{S+h}^* x_{S+h}) \\
&=&  \Re(z_{S+h}^* \wt\Delta_h x_{S+h}) +  \| z_{S+h} \|\,\| x_{S+h} \|  \ge 0.
\end{eqnarray*}
Combined with the extremality assumption, we obtain
$\Re(z^* \dot \Delta x) = 0$. This implies that $\wt\Delta_h$ has singular vectors $u_h$ and 
$v_h$, which completes the proof.
\end{proof}

\smallskip
\begin{remark} \rm
Theorem \ref{th:max2} allows us to restrict the perturbations in the structured spectral value set~\eqref{eq:svs} to those with rank-1 blocks, which was also shown in \cite{PD93}. Since the Frobenius and the matrix 2-norms of a rank-1 matrix are equal, we can equivalently search for extremizers within the submanifold  
\begin{eqnarray}
\BPCO & = & \bigl\{ \diag ( \delta_1 \Id_{r_1}, \ldots, \delta_S \Id_{r_S},
\Delta_1, \ldots,\Delta_F ): \bigr.
\nonumber
\\[1mm]
& & \bigl. \ \delta_i \in \C, | \delta_i | = 1,\ 
\Delta_j \in \C^{m_j\times m_j}, \| \Delta_j \|_\F = 1 \bigr\}.
\label{eq:optsetC}
\end{eqnarray}
\end{remark}

\subsection{A system of ODEs to compute extremal points of $\Lameps^{\BPC}(M)$} 
\label{sec:odes}


In order to compute a local maximizer for $|\lambda|$, 
with $\lambda \in \Lameps^{\BPC}(M)$, we will first construct 
a matrix valued function $\Delta(t)$, where $\Delta(t) \in \BPCO$, 
such that a largest eigenvalue $\lambda(t)$ of $\eps M \Delta(t)$ 
has maximal local increase. We then derive a system of ODEs satisfied by this choice of $\Delta(t)$.


%
\subsection*{Orthogonal projection onto $\BPC$}

In the following, we make use of the Frobenius inner product
$\langle A, B \rangle = {\rm trace}( A^* B)$
for two $m \times n$ matrices $A,B$. We let
\begin{equation}
C_\BPC = P_\BPC(C). 
\label{def:PB}
\end{equation}
denote the orthogonal projection, with respect to the Frobenius inner product, of a matrix $C \in \C^{n\times n}$ onto 
$\BPC$. To derive a compact formula for this projection, we use the pattern matrix 
\begin{equation}
\One_\BPC =  \diag \left( \One_{r_1}, \ldots, \One_{r_S},
\One_{m_1}, \ldots,\One_{m_F} \right),
\label{eq:OneBP}
\end{equation}
where $\One_{d}$ denotes the $d \times d$-matrix of all ones.
%
\begin{lemma} \label{lem:projS}
For $C \in \C^{n\times n}$, let 
\[
C \odot \One_\BPC = \diag \left( C_1, \ldots, 
C_{S+F} \right) 
\]
denote the block diagonal matrix obtained by entrywise multiplication of $C$ with
the matrix $\One_\BPC$ defined in~\eqref{eq:OneBP}.
Then the orthogonal projection of $C$ onto $\BPC$ is given by
\begin{equation}
C_\BPC = P_\BPC(C)  =   \diag \left( \gamma_1 \Id_{r_1}, \ldots, \gamma_S \Id_{r_S},
\Gamma_1, \ldots,\Gamma_F \right)
\label{eq:PBP}
\end{equation}
where $\gamma_i = \tr(C_i)/{r_i}, \ i=1,\ldots,S$,
and $\Gamma_1 = C_{S+1},\ldots,\Gamma_F = C_{S+F}$. 
\end{lemma} 
\begin{proof}
The result follows directly from the fact that
\[
\gamma_* = \arg\min\limits_{\gamma \in \C} \| E - \gamma \Id_r \|_\F = \frac{1}{r} 
\tr(E)
\] 
holds for every $E \in \C^{r\times r}$.
\end{proof}

If $C = u v^*$ is a rank-$1$ matrix, with the partitioning
\begin{eqnarray*}
u  =  \left( u_1^{\rm T} \ \ldots \ u^{\rm T}_S, \ u^{\rm T}_{S+1} \ \ldots \ u_{S+F}^{\rm T} \right)^{\rm T}, \quad 
v  =  \left( v_1^{\rm T} \ \ldots \ v^{\rm T}_S, \ v^{\rm T}_{S+1} \ \ldots \ v_{S+F}^{\rm T} \right)^{\rm T},
\end{eqnarray*} 
then the diagonal blocks $\Gamma_j = u_{S+j} v_{S+j}^*$ of the orthogonal projection are again rank-$1$ matrices
and, moreover, $\gamma_i = v_i^* u_i/r_i$. 


\subsection*{The local optimization problem}

Let us recall the setting from Section~\ref{sec:extrc}: We assume that $\lambda = |\lambda| \e^{\iu \theta}$
is a simple eigenvalue with eigenvectors $x,y$ normalized such that
\begin{eqnarray}
&& \| y \| = \| x \| = 1, \qquad y^* x = |y^* x| \e^{-\iu \theta}. 
\label{eq:scalyx}
\end{eqnarray}
As a consequence of Lemma \ref{lem:eigderiv}, see also~\eqref{eq:lambda2}, we have
\begin{eqnarray}
\frac d{dt}|\lambda|^2 
& = &  2 |\lambda| \Re\Bigl(\frac{z^*\dot \Delta x} {\e^{\iu \theta} y^* x} \Bigr)
= \frac{2 |\lambda|}{|y^* x|} \Re(z^*\dot \Delta x),
\label{eq:lambda3}
\end{eqnarray}
where $z = M^* y$ and the dependence on $t$ is intentionally omitted. 

Letting $\Delta\in \BPCO$, with $\BPCO$ as in~\eqref{eq:optsetC}, we now aim at determining a direction $\dot \Delta = \Z$ 
that locally maximizes the increase of the modulus of $\lambda$. This amounts to
%
determining
\begin{eqnarray}
\Z & = & \diag \left( \omega_1 \Id_{r_1}, \ldots, \omega_s \Id_{r_S},
\Omega_1,\, \ldots, \Omega_F \right)
\label{eq:Z}
\end{eqnarray} 
as a solution of the optimization problem
\begin{equation}
\begin{split}
\Z_* = & \arg\max \left\{ \Re( z^* \Z x ):\, \text{$\Z$ takes the form~\eqref{eq:Z}}\right\}   \\
\text{subject to }  & \quad \Re(\overline \delta_i \omega_i) = 0, \qquad \qquad \quad \ i=1,\ldots,S, \\
\text{ and } & \quad \Re \la \Delta_j,\Omega_j \ra = 0, \qquad \qquad j=1,\ldots,F.
\end{split}
\label{eq:propt} 
\end{equation}
The target function in~\eqref{eq:propt} follows from~\eqref{eq:lambda3}, while the constraints in~\eqref{eq:Z} and~\eqref{eq:propt} ensure that $\Z$ is in the tangent space of $\BPCO$ at $\Delta$. In particular,~\eqref{eq:propt} implies that the the norms of the blocks of $\Delta$ are conserved.
Note that~\eqref{eq:propt} only becomes well-posed after imposing an additional normalization on the norm of $Z$. The scaling chosen in the following lemma aims at 
$Z \in \BPCO$.

\begin{lemma}
\label{optlemmaR}
With the notation introduced above and $x,z$ partitioned as in~\eqref{eq:partitionxz}, a solution 
of the optimization problem 
\eqref{eq:propt} is given by 
\begin{eqnarray}
\Z_* & = & \diag \left( \omega_1 \Id_{r_1}, \ldots, \omega_S \Id_{r_S},
\Omega_1,\, \ldots, \Omega_F \right),
\nonumber
\end{eqnarray}
with
\begin{eqnarray}
\omega_i & = & 
\nu_i \left( x_i^* z_i - \Re\left( x_i^* z_i \overline{\delta}_i \right) \delta_i \right), 
\quad i=1,\ldots,S  
\label{eq:etaic}
\\[2mm]
\Omega_j & = & \zeta_j 
\left( z_{S+j} x_{S+j}^* - \Re \langle \Delta_j, z_{S+j} x_{S+j}^* \rangle \Delta_j \right),
\quad j=1,\ldots,F.
\label{eq:Omegaj}
\end{eqnarray}
Here, $\nu_i > 0$ is the reciprocal of the absolute value of the right-hand side in 
\eqref{eq:etaic}, if this is different from zero, and $\nu_i = 1$ otherwise. Similarly, $\zeta_j > 0$ 
is the reciprocal of the Frobenius norm of the matrix on the right hand side in \eqref{eq:Omegaj}, 
if this is different from zero, and $\zeta_j = 1$ otherwise. If all right-hand sides are different from zero
then $\Z_* \in \BPCO$.
\end{lemma}
\begin{proof}
The equality
\[
z^*\Z x = 
\sum\limits_{i=1}^{S} \omega_i z_i^* x_i + 
\sum\limits_{j=1}^{F} z_{S+j}^* \Omega_j x_{S+j} = \sum\limits_{i=1}^{S} \omega_i \langle z_i, x_i \rangle + 
\sum\limits_{j=1}^{F} \langle z_{S+j} x_{S+j}^*, \Omega_j \rangle
\]
implies that the maximization problem~\eqref{eq:propt} decouples, which allows us to maximize for each block of $Z$ individually.

For a full block $\Omega_j$, the term $\langle z_{S+j} x_{S+j}^*, \Omega_j \rangle$ is maximized by
the orthogonal projection of $z_{S+j} x_{S+j}^*$ onto the (real linear) subspace $\{ \Omega \in \C^{m_j\times m_j}:\, \Re\la \Delta_j,\Omega \ra = 0 \}$. This gives~\eqref{eq:Omegaj}, with the scaling chosen such that $\|\Omega_j\|_\F = 1$ unless $\Omega_j =0$. 

For a block $\omega_i I_{r_i}$, the term $\omega_i z_i^* x_i$ is maximized by projecting $x_i^* z_i$ onto 
$\{ \omega_i \in \C:\, \Re(\overline \delta_i \omega_i) = 0\}$. This gives~\eqref{eq:etaic}, with the scaling chosen such that $|\delta_i| = 1$ unless $\delta_i=0$.
\end{proof}
\smallskip
\begin{corollary} \label{cor:optlemmaR}
The result of Lemma {\rm \ref{optlemmaR}} can be expressed as
\begin{eqnarray}
\Z_{*} & = & D_1 P_\BPC \left( z x^* \right) - D_2 \Delta
\label{eq:Zopt}
\end{eqnarray}
where $P_\BPC(\cdot)$ is the orthogonal projection from Definition {\rm \ref{def:PB}},
and $D_1, D_2 \in \BPC$ are diagonal matrices with $D_1$ positive.
\end{corollary}
\begin{proof}
The statement is an immediate consequence of Lemma \ref{lem:projS}.
\end{proof}
\smallskip

\subsection*{The system of ODEs}

Lemma~\ref{optlemmaR} and Corollary~\ref{cor:optlemmaR} suggest to consider the following differential equation on the 
manifold $\BPCO$:
\begin{equation}\label{Fode}
\dot \Delta= D_1 P_\BPC ( z x^* ) - D_2 \Delta,
\end{equation} 
where $x(t)$ is an eigenvector, of unit norm, associated to a simple eigenvalue $\l(t)$ of 
$\eps M \Delta(t)$ for some fixed $\eps>0$. Note that $z(t),D_1(t),D_2(t)$ depend on $\Delta(t)$ as well.
The differential equation~\eqref{Fode} is a gradient system because, by definition, the right-hand side is the projected gradient of $Z\mapsto \Re(z^* Z x)$.

The following result follows directly from Lemmas~\ref{lem:eigderiv} and \ref{optlemmaR}.
\begin{theorem} \label{thm:monotone}
Let $\Delta(t) \in \BPCO$ satisfy the differential equation \eqref{Fode}.
If $\l(t)$ is a simple eigenvalue of $\eps M \Delta(t)$, then $|\l(t)|$ increases 
monotonically.
\end{theorem}

The following lemma establishes a useful property for the analysis of stationary points
of~\eqref{Fode}. 
\begin{lemma}
\label{lem:nonzero}
Let $\Delta(t) \in \BPCO$ satisfy the differential equation \eqref{Fode}.
If $\lambda(t)$ is a nonzero simple eigenvalue of $\eps M \Delta(t)$, with right and left 
eigenvectors $x(t)$ and $y(t)$ scaled according to \eqref{eq:scalyx}, then
\begin{equation}
P_\BPC\big(z(t) x(t)^* \big) \neq 0,
\label{eq:nonzero}
\end{equation}
where $z(t) = M^*y(t)$.
\end{lemma}
\begin{proof}
For convenience, we again omit the dependence on $t$ and let $\lambda = |\lambda| \e^{\iu \theta}$.
Assume -- in contradiction to the statement -- that  $P_\BPC \big( z x^* \big) = 0$. Because of the block diagonal structure of $\Delta$, this implies
\begin{equation}
 \Re \big\langle z x^*, \eps \Delta \big\rangle = 
   \Re \big\langle P_\BPC \left( z x^* \right), \eps \Delta \big\rangle = 0.
\label{eq:pzero}
\end{equation}
On the other hand,
\[
\Re \big\langle z x^*, \eps \Delta \big\rangle  = 
\Re \big\langle y x^*, \eps M \Delta \big\rangle =  
\Re \left( y^* \eps M \Delta x \right) = 
\Re \left( |\lambda| \e^{\iu \theta} y^*x \right).
\]
Exploiting the normalization \eqref{eq:scalyx} and the simplicity of $\lambda$, we obtain
$
\Re \big\langle z x^*, \eps \Delta \big\rangle  =  |\lambda|\, |y^* x| > 0.
$
This, however, contradicts \eqref{eq:pzero}.  
\end{proof}
\smallskip

The differential equation~\eqref{Fode} can be expressed in terms of the blocks of $\Delta$, that is, through 
$\{ \delta_i \}_{i=1}^{S}$ and $\{ \Delta_j \}_{j=1}^{F}$, as follows:
\begin{equation}\label{Fode2}
\begin{aligned}
\dot{\delta}_i & =  
\nu_i 
\Bigl( x_i^* z_i - \Re\left( x_i^* z_i \overline{\delta}_i \right) \delta_i \Bigr), & i = 1,\ldots,S \\
\dot{\Delta}_j & =  
\eta_j 
\Bigl( z_{S+j} x_{S+j}^* - \Re \langle \Delta_j, z_{S+j} x_{S+j}^* \rangle \Delta_j \Bigr), \ &
j = 1,\ldots,F
\end{aligned}
\end{equation}
with the scalars $\nu_i$, $\eta_j$ defined in Lemma~\ref{optlemmaR}.

Because of $|\delta_i| = 1$, we can reparametrize $\delta_i = \e^{\iu \beta_i}$ and rewrite the first set of equations in~\eqref{Fode2} 
as a system of ODEs in $\beta_i \in \R$. 
Setting $\gamma_i = \arg( z_i^* x_i )$, we obtain
\begin{eqnarray}
\iu \dot{\beta_i} \e^{\iu \beta_i} & = & \nu_i |z_i^* x_i| \left( \e^{-\iu \gamma_i} - 
\Re\left( \e^{-\iu (\gamma_i + \beta_i)} \right)  \e^{\iu \beta_i} \right),
\nonumber
\end{eqnarray}  
which gives
$
\dot{\beta_i}  =  -\nu_i |z_i^* x_i| \sin(\gamma_i +\beta_i).
$
With the normalization imposed by $\nu_i$, this finally yields
\[
\dot{\beta_i}  =  -{\rm sign} \left( \sin(\gamma_i + \beta_i) \right).
\]
In particular, this means that $\dot{\beta_i} = 0$ if and only if $\gamma_i + \beta_i = 0, \pm \pi$;
maximizers correspond to $\beta_i = -\gamma_i$.

\begin{remark}
The choice of $\nu_i$, $\eta_j$ originating from Lemma~\ref{optlemmaR}, to achieve unit norm of all blocks in~\eqref{eq:Zopt}, is completely arbitrary.
Other choices would be also acceptable and investigating an optimal one
in terms of speed of convergence to stationary points would be an interesting issue.
\end{remark}

The following result characterizes stationary points of \eqref{Fode}.
\begin{theorem}
\label{stat:RF}
Assume that $\Delta(t)$ is a solution of \eqref{Fode} and $\l(t)$ is a largest simple nonzero eigenvalue of 
$\eps M \Delta(t)$ with right/left eigenvectors $x(t)$, $y(t)$. Moreover, suppose that
Assumptions~\eqref{eq:A1} and~\eqref{eq:A2} hold for $x(t)$ and 
$z(t)=M^* y(t)$. 
Then 
\begin{equation} \label{eq:equivalence}
\frac{d}{d t} |\l(t)|^2=0 \ \Longleftrightarrow \ 
\dot  \Delta(t)=0  \ \Longleftrightarrow \ 
\Delta(t) = D P_\BPC \left(z(t) x(t)^* \right),
\end{equation}
for a specific real diagonal matrix $D \in \BPC$. 
Moreover if $\l(t)$ has (locally) maximal modulus over the set $\Lambda_\eps^{\BPC}(M)$ then 
$D$ is positive. 


\end{theorem}
\begin{proof}
By~\eqref{eq:lambda3}, $\frac{d}{d t} |\l(t)|^2=0$ implies $\Re(z^* \dot  \Delta x) = 0$.
Inserting~\eqref{Fode2} shows that each block of $\dot  \Delta$ is necessarily zero and hence $\dot  \Delta = 0$.
The other direction of the first equivalence in~\eqref{eq:equivalence} is trivial.
The second equivalence in~\eqref{eq:equivalence} follows directly from~\eqref{Fode}. 
By Theorem~\ref{th:max}, all blocks of $\Delta(t)$ have norm $1$ and hence none of the scalars defining $D$ can be zero.
Thus, $D$ is nonsingular.

Assumping that $\l(t)$ has (locally) maximal modulus, we now prove positivity of $D$ by contradiction.
Suppose that the $(S+j)$th full block of $D$ is equal to $-\gamma_j\,\Id_{m_j}$ with $\gamma_j > 0$, implying  $\Delta_j = -\gamma_j z_{S+j} x_{S+j}^*$. Consider an ODE with the $(S+j)$th block 
$\dot \Delta_j = z_{S+j} x_{S+j}^*$ and initial datum $\Delta_j(0) = -\gamma_j z_{S+j} x_{S+j}^*$, 
while leaving all other blocks of $\Delta$ unaltered.
Such an ODE clearly decreases the norm of $\Delta_j$ for $t\le \bar t$, for some $\bar t > 0$ (implying that 
$\|\Delta(t)\|_\F$ does not exceed $1$). 
By the usual derivative formula from Lemma~\ref{lem:eigderiv} the largest eigenvalue 
$\lambda(t)$ of $\eps \Delta(t) M$ is such that $\frac{d}{d t} |\lambda(t)| > 0$, which contradicts 
local maximality. 

Similarly consider a repeated scalar block and assume that the $i$th block of $D$ is equal to 
$-\gamma_i\,\Id_{r_i}$ with $\gamma_i > 0$, which means $\delta_i = \gamma_i z_i^* x_i$. Consider, similarly
to previous case, an ODE with the $i$th block $\dot \delta_i = x_i^* z_i$ and initial datum
$\delta_i = -\gamma_i x_i^* z_i$. Again, $|\delta_i(t)|$ decreases in a sufficiently small time-interval
$[0,\bar t]$ and $|\lambda(t)|$ increases in the same interval, yielding again a contradiction.    
\end{proof}

\subsection{Projection of full blocks on rank-$1$ manifolds}
\label{sec:proj}

In order to exploit the rank-$1$ property of extremizers established in Theorem~\ref{th:max2}, we can proceed in complete 
analogy to \cite{GL11} in order to obtain for each full block an ODE on the manifold 
$\cM_1$ of (complex) rank-1 matrices. 
We express $\Delta_j\in\cM_1 \subset \C^{m_j\times m_j}$  as
$$
\Delta_j=\sigma_j p_j q_j^*, \quad
\dot \Delta_j = \dot\sigma_j p_j q_j ^* + \sigma_j \dot p_j q_j ^* + \sigma_j p_j \dot q_j ^*
$$
where $\sigma_j\in \C$ and $p_j,q_j\in \C^{m_j}$ have unit norm. The parameters
$\dot\sigma_j\in\C$, $\dot p_j,\dot q_j\in \C^{m_j}$ are uniquely determined 
by $\sigma_j,p_j,q_j$ and $\dot \Delta_j$ when imposing the orthogonality conditions 
$p_j ^*\dot p_j=0, \  q_j ^*\dot q_j=0$.

In the differential equation (\ref{Fode2}) we replace the right-hand side by its orthogonal 
projection onto the tangent space $T_{\Delta_j} \cM_1$ (and also remove the normalization constant) to obtain
\begin{equation}\label{odeDeltaj-1}
\dot \Delta_j = 
P_{\Delta_j} \left( z_{S+j} x_{S+j}^* - \Re \langle \Delta_j, z_{S+j} x_{S+j}^* \rangle \Delta_j \right).
\end{equation}
Note that the orthogonal projection of a matrix $Z\in\C^{m_j\times m_j}$ onto $T_{\Delta_j} \cM_1$ at 
$\Delta_j=\sigma_j p_j q_j ^* \in\cM_1$ is given by
\[
P_{\Delta_j}(Z) = Z - (\Id-p_j p_j ^*) Z (\Id-q_j q_j ^*).
\]

Following the arguments of \cite{GL11}, the equation $\dot \Delta_j=P_{\Delta_j}(Z)$ is equivalent to
\begin{eqnarray*}
\dot \sigma_j &=& p_j ^*  Z q_j
\nonumber\\
\dot p_j &=& (\Id-p_j p_j ^*) Z q_j \sigma_j^{-1}
\label{odes-C-1} \\
\dot q_j &=& (\Id-q_j q_j ^*) Z ^*  p_j {\overline\sigma_j}^{-1}.
\nonumber
\end{eqnarray*}
Inserting $Z= z_{S+j} x_{S+j}^* - \Re \langle \Delta_j, z_{S+j} x_{S+j}^* \rangle \Delta_j$, we obtain that 
the differential equation~\eqref{odeDeltaj-1} is equivalent to the following system of differential equations 
for $\sigma_j,p_j$ and $q_j$, where we set 
$\alpha_j=p_j ^* z_{S+j}\in \C$, $\beta_j=q_j ^* x_{S+j} \in \C$:
\begin{eqnarray}
\dot \sigma_j &=& \alpha_j\overline\beta_j - \Re(\overline\alpha_j \beta_j \sigma_j) \sigma_j
               =  \iu\,\Im(\alpha_j\overline\beta_j\overline\sigma_j)\sigma_j
\nonumber
\\
\dot p_j &=& ( z_{S+j}-\alpha_j p_j)\overline\beta_j \sigma_j^{-1}
\label{ode-spq-1} 
\\
\dot q_j &=& (x_{S+j}-\beta_j q_j){\overline\alpha_j} \,{\overline\sigma_j}^{-1}.
\nonumber
\end{eqnarray}
The derivation of this system of ODEs is straightforward; we refer the reader to \cite{GL13}
for details.

The monotonicity and the characterization of stationary points follows analogously
to those obtained for \eqref{Fode2}; we refer to~\cite{GL11} for the proofs.
As a consequence we can use the ODE~\eqref{ode-spq-1} instead of \eqref{Fode2} and gain in terms of computational complexity.

\subsection{An illustrative example}
\label{sec:illc}

Consider the matrix
\begin{eqnarray*}
M & = & \footnotesize \left( \begin{array}{rrrrr} 
-0.10 - 0.55\,\iu & -0.57 - 1.59\,\iu & -1.34 - 1.70\,\iu &  0.04 + 0.49\,\iu & -0.18 + 0.19\,\iu \\
-1.48 - 2.17\,\iu &  0.58 + 1.17\,\iu &  0.05 + 0.53\,\iu &  0.11 - 0.42\,\iu &  0.26 + 1.19\,\iu \\
-0.53 + 0.59\,\iu &  0.78 - 1.48\,\iu &  0.15             & -0.25 + 1.46\,\iu &  0.33 + 1.32\,\iu \\
 0.24 + 0.79\,\iu & -0.12 - 0.65\,\iu &  1.79 - 0.09\,\iu & -0.63 + 1.39\,\iu & -0.88 + 0.10\,\iu \\
-2.03 + 1.33\,\iu & -1.22 - 0.22\,\iu &  0.45 - 1.49\,\iu &  0.94 - 0.13\,\iu & -1.02 + 2.33\,\iu 
	\end{array} \right)
\end{eqnarray*}
and a perturbation set given by
\[
\BPC = \big\{ \diag \left( \delta_1 \Id_{1}, \delta_2 \Id_{1}, \Delta_1, \delta_3 \Id_{1}\right), \delta_1, \delta_2, \delta_3 \in \C, \Delta_1 \in \C^{2,2} \big\}.
\] 

Applying {\sc Matlab}'s {\tt mussv}, we obtain the perturbation $\widehat\eps \widehat\Delta$
(with $\| \widehat\Delta\|_2 = 1$)  
\begin{eqnarray*}
&& 
\widehat\Delta =
\left( \begin{array}{rrrr}
 \e^{-\iu\,0.48650737} & 0 & 0 & 0 \\                     
  0 & \e^{-\iu\,1.49644308}  & 0 & 0 \\                      
  0 & 0 & u_1 v_1^* & 0 \\
  0 & 0 & 0 & \e^{-\iu\,2.155849308} 
\end{array} \right) \\[2mm]
&&
u_1 = \footnotesize \left( \begin{array}{r} 
  0.41899793 + 0.68039781\,\iu \\
  0.06834008 - 0.59735180\,\iu 
\end{array} \right), \quad
v_1 = \left( \begin{array}{r} 
  0.52696073 \\
  0.70477030 + 0.47498548\,\iu 
\end{array} \right),
\end{eqnarray*}
and $\widehat\eps = 0.228726413$, which determines the following lower bound for the $\mu$-value:
$\mu_\BPC(M) \ge \lbo = 4.372035505$. 
  
Applying Algorithm~\ref{algo-ode} below we find the locally extremal perturbation
\begin{eqnarray*}
&& \Delta^\star =
\left( \begin{array}{rrrr}
 \e^{-\iu\,2.49033999} & 0 & 0 & 0 \\                     
  0 & \e^{\iu\,1.24640446}  & 0 & 0 \\                      
  0 & 0 & u_1 v_1^* & 0 \\
  0 & 0 & 0 & \e^{-\iu\,1.72494213} 
\end{array} \right) \\
&&
u_1 = \scriptsize \left( \begin{array}{r} 
  0.15703326 + 0.85130227\,\iu \\
  0.29626531 - 0.40354908\,\iu  
\end{array} \right), \quad
v_1 = \left( \begin{array}{r} 
  0.68793173 \\
  0.28357426 + 0.66808351\,\iu 
\end{array} \right).
\end{eqnarray*} 
and $\eps^\star = 0.222994978$, which determines the following lower bound for the $\mu$-value:
$\mu_\BPC(M) \ge \lbo = 4.484405922$. Thus, the lower bound has improved, in particular when taking into account that the upper bound computed by \mussv{} is
4.48638$\ldots$.

\begin{figure}[ht]
\centerline{
\includegraphics[scale=0.5,trim= 0mm 0.01mm 0mm 0mm]{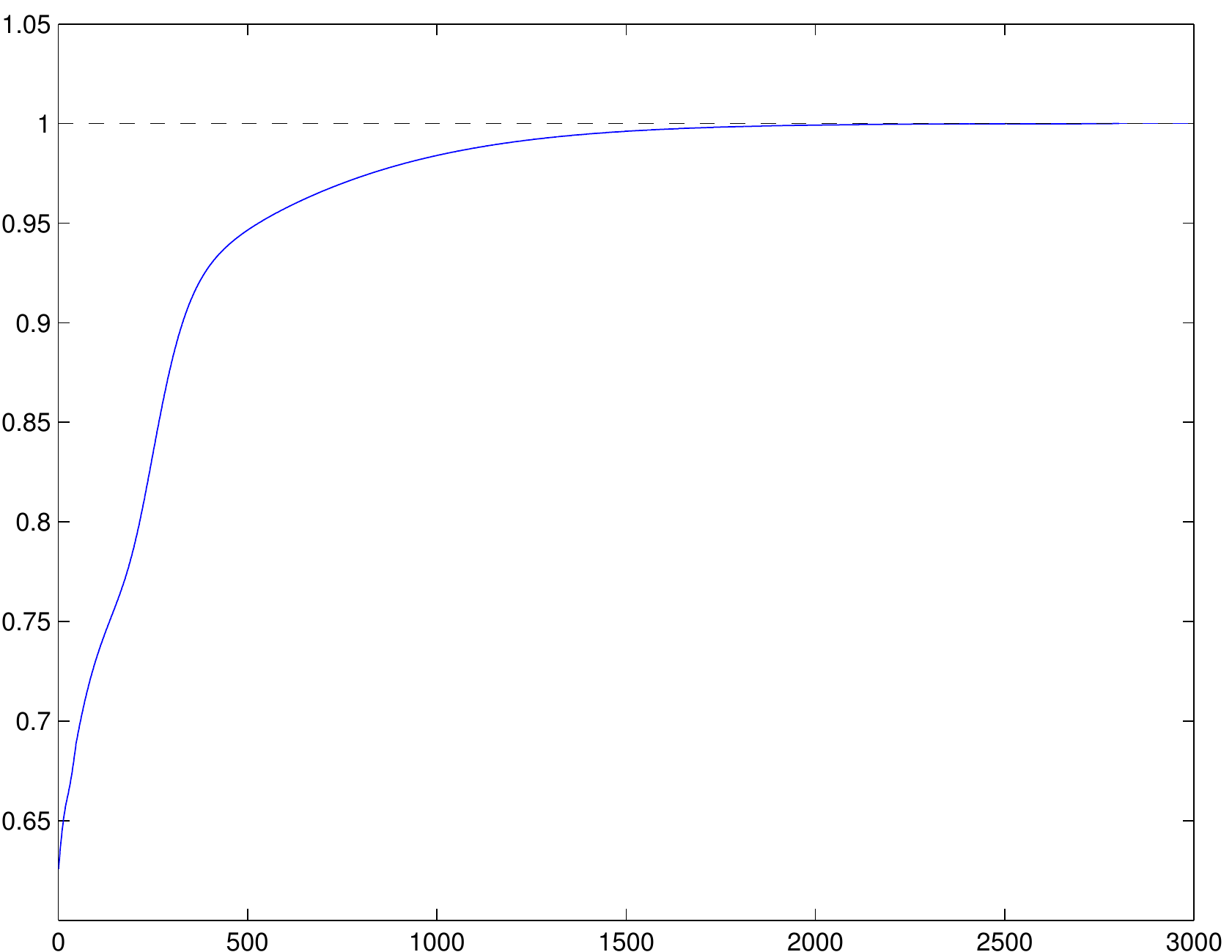}
}
\caption{Monotonic behavior of $\rho(\eps^\star M \Delta(t))$ along the
solution of the ODE \eqref{Fode} for the example from Section~\ref{sec:illc}. \label{fig1}}
\end{figure}

The behavior of the spectral radius of the matrix $\eps^\star M \Delta(t)$ along the
solution of the ODE is illustrated in Figure \ref{fig1}, which shows the monotonically
increasing behavior of $|\l(\eps)|$.

\section{General complex/real perturbations}
\label{sec:mix}

We now consider the more general case addressed by~\cite{PD93}, that is,
\begin{equation}
\label{eq:BPcr}
\BPR = \big\{ \diag \left( \delta_1 \Id_{r_1}, \ldots, \delta_s \Id_{r_S},
\Delta_1, \ldots,\Delta_F \right):\, \delta_i \in \C \mbox{(or $\R$)},
\Delta_j \in \C^{m_j \times m_j} \big\},
\end{equation}
where $\delta_i$ is either a complex or a real scalar.
Without loss of generality, we assume that the first
$S'$ repeated scalar blocks are complex while all other repeated scalar blocks are real.
Moreover, we set
\[
\BPR_1 \!=\! \big\{ \diag \left( \delta_1 \Id_{r_1}, \ldots, \delta_S \Id_{r_S},
\Delta_1, \ldots,\Delta_F \right) \in \BPR:\, |\delta_i|=1, \| \Delta_j \|_\F = 1 \big\}.
\]
This case differs qualitatively  from the purely complex case discussed in Section~\ref{sec:compl}, since it has to be formulated as a structured 
distance to singularity of the matrix $\Id - \eps M \Delta$. Due to the realness constraint
for some of the repeated scalar blocks, $\Delta \in \BPR$ does not
imply that ${\rm e}^{\iu \theta} \Delta \in \BPR$ for all $\theta \in [0, 2 \pi)$, which
means that the spectral value set $\Lambda_\eps(M)$ is generally not a disk. In turn, we need to address the minimization problem~\eqref{eq:cr} instead of the maximization problem~\eqref{eq:crcomplex}.

\subsection{Extremizers}
\label{sec:extrmix}

\begin{definition}
A matrix $\Delta \in \BPR$, such that $\| \Delta \|_2 \le 1$ and $\Id-\eps M \Delta$ has
a smallest eigenvalue that locally minimizes the modulus of $\Sigma_\eps^{\BP}(M)$
is called a {\em local extremizer}.
\end{definition}

We have the following result concerning local extremizers for the smallest (in modulus) 
complex number in $\Sigma_\eps^{\BP}(M)$.
\begin{theorem} \label{th:maxmix}
Let 
$$
\wt\Deltaopt = \diag \big( \wt\delta_1 \Id_{r_1}, \ldots, \wt\delta_{S'} \Id_{r_{S'}},
\wt\delta_{S'+1} \Id_{r_{S'+1}}, \ldots, \wt\delta_{S} \Id_{r_{S}},
\wt\Delta_1, \ldots,\wt\Delta_F \big), \quad \| \wt\Deltaopt \|_2 = 1, 
$$
be a local extremizer of $\Sigma_\eps^{\BP}(M)$.
Let $\lz  = |\lz| \e^{\imagunit \theta}$ be a simple smallest eigenvalue of the matrix $\Id-\eps M \wt\Deltaopt$, with 
the right and left eigenvectors $x$ and $y$ scaled such that $s = \e^{\imagunit \theta} y^* x > 0$.
Partitioning $x$ and $z = M^* y$ as in~\eqref{eq:partitionxz}, we assume that
\begin{eqnarray}
&& z_k^* x_k \neq 0 \quad \forall \ k=1,\ldots,S' 
\label{eq:M1}
\\
&& \Re(z_k^* x_k) \neq 0 \quad \forall \ k=S'+1,\ldots,S 
\label{eq:M2}
\\
&& \|z_{S+h}\|_2 \cdot \|x_{S+h}\|_2 \ne 0 \quad \forall \ h=1,\ldots,F
\label{eq:M3}
\end{eqnarray}
hold.
Then
$$
| \wt\delta_{k} | = 1 \quad \forall \ k=1,\ldots,S
\quad \mbox{and} \quad \| \wt{\Delta}_{h} \|_2 = 1 \quad \forall h=1,\ldots,F,
$$
that is, all blocks of $\wt\Deltaopt$ have unit $2$-norm.
\end{theorem}
\begin{proof}
\rm 
The proof is analogous to the proof of Theorem \ref{th:max}. The only substantial difference is caused by repeated real scalar blocks. 
To address this case, suppose that $| \wt\delta_{k} | < 1$ for some $S'+1 \le k \le S$ with $\delta_k \in \R$.
Let us consider the matrix valued function
\[
\Delta(t) = \diag \left( \wt\delta_{1} \Id_1,\ldots,\wt\delta_{k} \Id_k - 
t\, \Re(x_k^*z_k) \Id_{k}, 
\ldots, \delta_{S} \Id_S, \wt\Delta_1, \ldots,\wt\Delta_F \right),
\]
which satisfies $\Delta(0)=\wt\Deltaopt$ and $\| \Delta(t) \|_2 \le 1$ for $t$ sufficiently small. 
Making use of Lemma~\ref{lem:eigderiv} and~\eqref{eq:M2}, we obtain
\[
\frac{d}{d t} |\lz(t)|^2 \Big|_{t=0} =  
-2 \eps \frac{|\lz|}{s} \Re(z_k^* x_k)^2 < 0,
\]
which contradicts the extremality of $\wt\Deltaopt$ and thus completes the proof.
\end{proof}

\subsection{A system of ODEs to compute extremal points of $\Sigma_\eps^{\BP}(M)$}
\label{sec:gsmix}

We can apply a procedure similar to the one developed in the Section~\ref{sec:odes} 
to develop a system of ODEs for solving~\eqref{eq:cr}, that is, $\arg\min\limits_{\z \in \Sigma_\eps^{\BP}(M)} |\z|$.

The first step is to determine a steepest descent direction for $|\lz|$, which then yields the right-hand side of a suitable gradient system
for computing a local minimizer of the modulus of $\Sigma_\eps(M)$.
%
For this purpose, we use the same normalization~\eqref{eq:scalyx}
for the
eigenvectors $x$ and $y$ associated to a simple smallest eigenvalue $\lz$ of $\Id - \eps M \Delta$, 
and recall from~\eqref{eq:lambda3} that
\begin{eqnarray}
\frac d{dt}|\lz|^2 
& = &  -
\frac{2 |\lz|}{|y^* x|} \Re(z^*\dot \Delta x), \qquad \mbox{with} \ z = M^* y.
\nonumber
\label{eq:lambdacr3}
\end{eqnarray}
%
Rewriting the (constrained) minimization of this expression in terms of the blocks of
\begin{equation}\label{eq:Z2}
\dot\Delta = \Z  =  \diag \left( \omega_1 \Id_{r_1}, \ldots, \omega_s \Id_{r_S},
\Omega_1,\, \ldots, \Omega_F \right) \in \BPR
\end{equation}
yields the following optimization problem: 
\begin{equation}\label{eq:proptcr}
\begin{aligned}
\Z_* = & \arg\max \left\{ \Re( z^* \Z x ):\, \text{$\Z$ takes the form~\eqref{eq:Z2}}\right\}   \\
\text{subject to }  & \quad \Re(\delta_i^* \omega_i) = 0, \qquad\qquad\quad\ i=1,\ldots,S', \\
\text{ and } & \quad \Re \la \Delta_j,\Omega_j \ra = 0, \qquad \qquad j=1,\ldots,F, \\
\text{ and } & \quad \delta_\ell\,\omega_\ell \le 0 \quad \mbox{if} \ \delta_\ell=\pm 1, \qquad 
\!\! \ell=S'+1,\ldots,S. 
\end{aligned}
\end{equation}
As before, the first two constraints imply the conservation of the Frobenius norms for the full and
repeated complex scalar blocks of $\Delta$. The third constraint prevents $|\delta_\ell|$ from
exceeding $1$ for repeated real scalar blocks.
To make~\eqref{eq:proptcr} well-posed, we need to impose a normalization on $Z$ and the following lemma aims at $Z \in \BPR_1$, whenever this is possible.
\begin{lemma}
\label{optlemmaRcr}
With the notation introduced above and $x,z$ partitioned as in~\eqref{eq:partitionxz}, a solution of the 
optimization problem~\eqref{eq:proptcr} is given by 
\begin{eqnarray}
\Z_* & = & \diag \left( \omega_1 \Id_{r_1}, \ldots, \omega_s \Id_{r_S},
\Omega_1,\, \ldots, \Omega_F \right)
\nonumber
\end{eqnarray}
with
\begin{eqnarray}
\omega_i & = & 
\nu_i \left( x_i^* z_i - \Re\left( x_i^* z_i \overline{\delta}_i \right) \delta_i \right), 
\quad i=1,\ldots,S'  
\label{eq:etaicc} \\
\omega_\ell & = & 
\left\{ \begin{array}{rl} 
1 & \mbox{\rm if} \ \Re(z_\ell^* x_\ell) > 0 \ \mbox{\rm and} \ \delta_\ell > -1 \\
 -1 & \mbox{\rm if} \ \Re(z_\ell^* x_\ell) < 0 \ \mbox{\rm and} \ \delta_\ell <  1, \\
 0 & \mbox{\rm otherwise}
\end{array}
\right. \quad \ell=S'+1,\ldots,S
\label{eq:etair}
\\
\Omega_j & = & \eta_j 
\left( z_{S+j} x_{S+j}^* - \Re \langle \Delta_j, z_{S+j} x_{S+j}^* \rangle \Delta_j \right),
\quad j=1,\ldots,F.
\label{eq:Omegacrj}
\end{eqnarray}
Here, $\nu_i > 0$ is the reciprocal of the absolute value of the right-hand side in 
\eqref{eq:etaicc}, if this is different from zero, and $\nu_i = 1$ otherwise; $\eta_j > 0$ 
is the reciprocal of the Frobenius norm of the matrix on the right hand side in \eqref{eq:Omegacrj}, 
if this is different from zero, and $\eta_j = 1$ otherwise. 
\end{lemma}
\begin{proof}
The equality 
\begin{eqnarray*}
z^*\Z x&= & 
\sum\limits_{i=1}^{S'} \omega_i z_i^* x_i + 
\sum\limits_{\ell=S'+1}^{S} \omega_\ell z_\ell^* x_\ell +
\sum\limits_{j=1}^{F} z_{S+j}^* \Omega_j x_{S+j} 
 \\
&= &
\sum\limits_{i=1}^{S'} \omega_i \langle z_i, x_i \rangle + 
\sum\limits_{\ell=S'+1}^{S} \omega_\ell \langle z_i, x_i \rangle +
\sum\limits_{j=1}^{F} \langle z_{S+j} x_{S+j}^*, \Omega_j \rangle. 
\end{eqnarray*}
allows us to maximize for each block of $Z$ individually. The expressions~\eqref{eq:etaicc} and~\eqref{eq:Omegacrj} follow from the proof of Lemma~\ref{optlemmaR}.
The expression~\eqref{eq:etair} trivially maximizes the real part of $\omega_\ell z_\ell^* x_\ell$ among $\omega_\ell = \pm 1$.
\end{proof}

\subsection{The system of ODEs}
\label{sec:odesmix}

We use the result of Lemma~\ref{optlemmaRcr} to build a gradient system for $\Delta$ to find a local minimizer of $|\lz|$.
In terms of the blocks of $\Delta$, we obtain the following system of differential equations:
%
\begin{equation} \label{Fodecr2}
 \begin{aligned}
\dot{\delta}_i & =  \nu_i \Bigl( x_i^* z_i - \Re\left( x_i^* z_i \overline{\delta}_i \right) \delta_i \Bigr), & i = 1,\ldots,S'
\\
\dot{\delta}_\ell & =  {\rm sign} \Bigl( \Re\left(z_\ell^* x_\ell \right) \Bigr)
\One_{(-1,1)}(\delta_\ell), & \ell = S'+1,\ldots,S
\\
\dot{\Delta}_j & =  \eta_j 
\Bigl( z_{S+j} x_{S+j}^* - \Re \langle \Delta_j, z_{S+j} x_{S+j}^* \rangle \Delta_j \Bigr),
& j = 1,\ldots,F,
 \end{aligned}
\end{equation}
where $\delta_i \in \C$ for $i=1,\ldots,S'$, $\delta_\ell \in \R$ for $\ell=S'+1,\ldots,S$, and
$\One_{E}(\cdot)$ is the characteristic function for a set $E$.

Expressing $\delta_i = \e^{\iu \beta_i}$ the first equation in~\eqref{Fodecr2} can again be rewritten as 
\[
\dot{\beta_i}  = -{\rm sign} \left( \sin(\gamma_i + \beta_i) \right),
\]
which means that $\dot{\beta_i} = 0$ if and only if $\gamma_i + \beta_i = 0, \pm \pi$;
extremizers correspond to $\beta_i = -\gamma_i$.

A system of ODEs that exploits the rank-1 property of the full blocks in extremizers can be derived in a fashion completely analogous to Section~\ref{sec:proj}.

\subsection{An illustrative example}
\label{sec:illcr}

Consider the matrix
\begin{eqnarray*}
M & = & \scriptsize \left( \begin{array}{rrrrr}
  -1.54 - 1.28 \iu &  -0.56 + 0.57 \iu &  -0.03 - 0.63 \iu & -0.64 - 0.55 \iu &   0.46 - 0.22 \iu \\
  -1.08 + 1.91 \iu &   1.16 - 0.08 \iu &  -0.41 - 0.13 \iu &  0.04 - 0.06 \iu &  -0.01 - 0.71 \iu \\
   0.11 - 2.16 \iu &   0.53 + 0.79 \iu &  -0.33 + 0.26 \iu &  0.44 + 0.02 \iu &   0.20 + 0.96 \iu \\
   0.52 + 0.29 \iu &   2.38 + 0.09 \iu &  -0.03 + 0.06 \iu &  0.01 + 1.12 \iu &   0.51 - 0.77 \iu \\
  -1.30 + 0.34 \iu &  -1.72 + 0.14 \iu &   1.02 + 1.34 \iu &  0.35 - 0.75 \iu &   0.48 + 0.04 \iu 
	\end{array} \right)
\end{eqnarray*}
and a perturbation set given by
\[
\BPR = \big\{ \diag \left( \delta_1 \Id_{1}, \delta_2 \Id_{1}, \delta_3 \Id_{1}, \delta_4 \Id_{2} \right), \delta_1, \delta_2 \in \R, \ \delta_3, \delta_4 \in \C \big\}.
\] 
Applying the {\sc Matlab}'s \mussv, we obtain the perturbation $\widehat\eps \widehat\Delta$
with  
\[
\widehat\Delta = \scriptsize 
\left( \begin{array}{rrrrr}
 -1 & 0 & 0 & 0 & 0 \\                     
  0 & 1 & 0 & 0 & 0 \\                      
  0 & 0 & \e^{-\iu\,0.91357833} & 0 & 0 \\
  0 & 0 & 0 & \e^{-\iu\,2.076961991} & 0 \\
  0 & 0 & 0 & 0 & \e^{-\iu\,2.076961991}
\end{array} \right).
\]
and $\widehat\eps = 0.30300829$, which yields the lower bound 
$\mu_\BPR(M) \ge \lbo = 3.300239739$. 
Applying the algorithms presented in this article we find the same solution $\lbn=\lbo$. 

\begin{figure}[ht]
\centerline{
\includegraphics[scale=0.5,trim= 0mm 0.01mm 0mm 0mm]{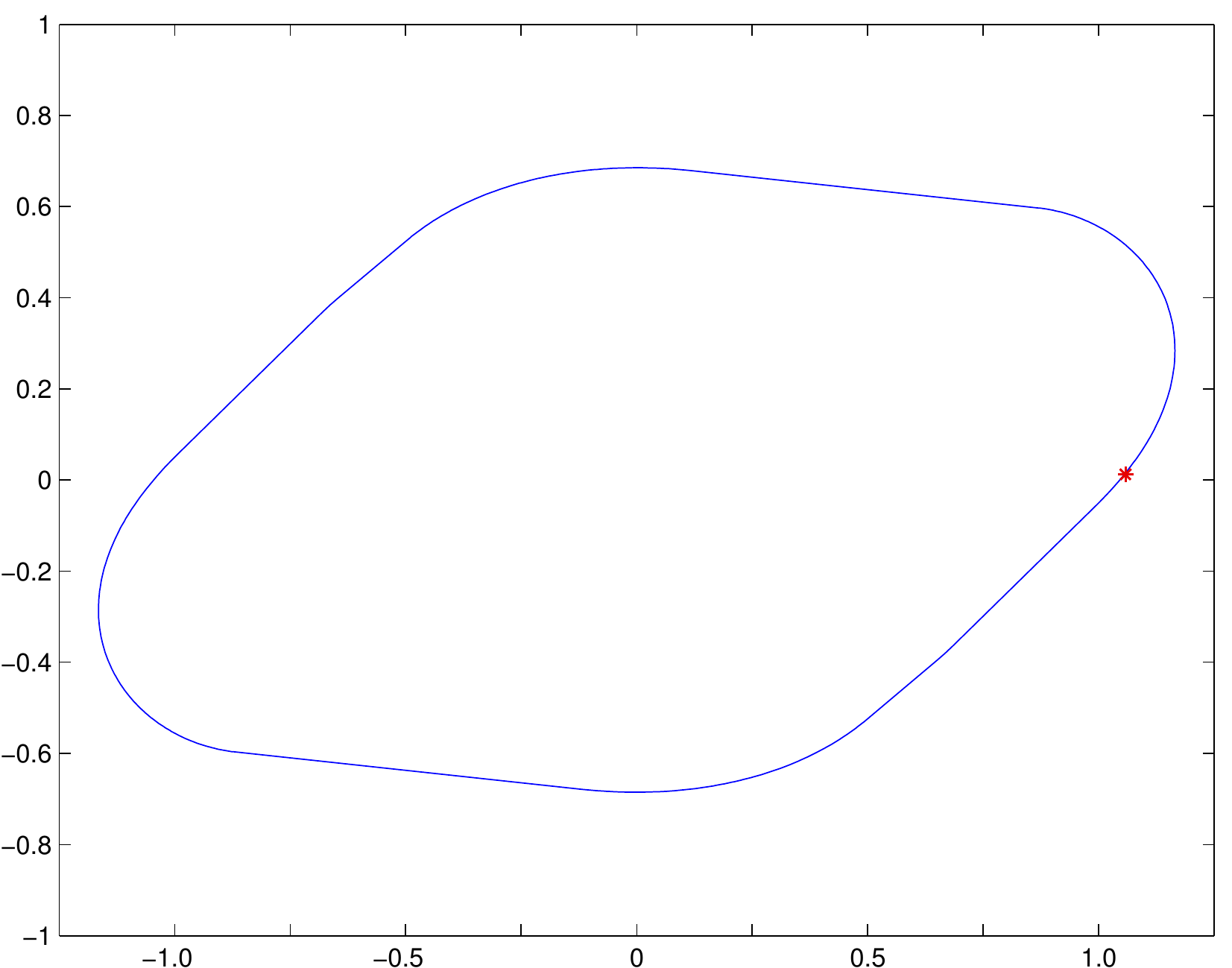}
}
\caption{Boundary of the set $\Lambda_\eps^{\BPR}(M)$ for $\eps= 0.30300829$. The point
$\z=1$ (the red asterisk in the picture) correctly lies on the boundary of $\Lambda_\eps^{\BPR}(M)$. \label{fig2}}
\end{figure}

Intensively sampling the set of all possible perturbations indicates that the computed value $\lbo$ yields the exact value of $\mu_\BPR(M)$; see also Figure~\ref{fig2}.

\subsection{Choice of initial value matrix and $\eps_0$}
\label{sec:delta0}

In our two-level algorithm for determining $\eps$, we use the perturbation $\Delta$ obtained for the previous value $\eps$ as the initial value matrix for the system of ODEs~\eqref{Fodecr2}. However, it remains to discuss a suitable choice of the initial values $\Delta(0) = \Delta_0$ and $\eps_0$ in the very beginning of the algorithm.

For the moment, let us assume that $M$ is invertible and write $$\Id - \eps_0 M \Delta_0 = M ( M^{-1} - \eps_0 \Delta_0 ),$$ 
which we aim to have as close as possible to singularity. To determine $\Delta_0$, we perform an asymptotic analysis around $\eps_0 \approx 0$.
For this purpose, let us consider the matrix valued function
\[
G(\tau) = M^{-1} - \tau \Delta_0,
\]
and let denote $\chi(\tau)$ denote an eigenvalue of $G(\tau)$ with smallest modulus.
Letting $x$ and $y$ denote the right and left eigenvectors corresponding to $\chi(0) = \chi_0 = |\chi_0| \e^{\iu \theta}$, 
scaled such that $\e^{\imagunit \theta} y^*x > 0$, Lemma~\ref{lem:eigderiv} implies
\begin{eqnarray*}
 \frac{d}{d \tau} |\chi(\tau)|^2 \Big|_{\tau=0}  &=& 
2\,\Re(\overline\chi \dot\chi) =
-2\,\Re\Bigl(\overline\chi\, \frac{y^* \Delta_0 x} {y^*x} \Bigr) \\
 &=&  
-2 |\chi_0| \Re\Bigl(\frac{y^* \Delta_0 x} {\e^{\iu \theta} y^*x} \Bigr) =
- \frac{2 |\chi_0|}{|y^* x|} \Re\langle y x^*, \Delta_0 \rangle.
\nonumber
\end{eqnarray*}
In order to have the locally maximal decrease of $|\chi(\tau)|^2$ at $\tau=0$ we choose 
\begin{eqnarray}
\Delta_0 = D\,P_\BPR(y x^*), 
\label{eq:Delta0}
\end{eqnarray} 
where the positive diagonal matrix $D$ is chosen such that $\Delta_0 \in \BPR_1$. 
This is always possible under the genericity assumptions~\eqref{eq:M1}--\eqref{eq:M3}.
The orthogonal projector $P_\BPR$ onto $\BPR$  can be expressed in analogy to~\eqref{eq:PBP} for $P_\BPC$, with the notable difference that
$\gamma_\ell = \Re(\tr(C_\ell))/r_\ell$ for $\ell=S'+1,\ldots,S$.
Note that there is no need to form $M^{-1}$; $x$ and $y$ can be obtained as the eigenvectors associated to a largest
eigenvalue of $M$. However, attention needs to be paid to the scaling. Since the largest eigenvalue of $M$ is
$\frac{1}{|\chi_0|} \e^{-\imagunit \theta}$, $y$ and $x$ have to be scaled accordingly. 

A possible choice for $\eps_0$ is obtained by solving the following simple linear equation, resulting from the first order expansion of the eigenvalue at $\tau=0$:
\begin{equation*}
|\chi(\eps_0)|^2 + \frac{d}{d \tau} |\chi(\tau)|^2 \Big|_{\tau=0} \eps_0  = 0. 
\end{equation*}
This gives
\begin{equation}
\eps_0 = \frac{|\chi_0|\, |y^* x|}{2 \Re\langle y x^*, \Delta_0 \rangle} = \frac{|\chi_0|\, |y^* x|}{2 \| P_\BP(y x^*) \|}.
\label{eq:choice1}
\end{equation} 
This can be improved in a simple way by computing this expression for $\eps_0$ for several eigenvalues of $M$ (say, the
$m$ largest ones) and taking the smallest computed $\eps_0$. For a sparse matrix $M$, the \matlab
function \verb@eigs@ (an interface for ARPACK, which implements the implicitly restarted  Arnoldi Method 
\cite{LS96,LSY98}) allows to efficiently compute a predefined number $m$ of Ritz values.

Another possible, very natural choice for $\eps_0$ is given by
\begin{equation}
\eps_0 = \frac{1}{\overline{\mu}_\BPR(M)}
\label{eq:choice2}
\end{equation}
where $\overline{\mu}_\BPR(M)$ is the upper bound for the SSV computed by the \matlab function \mussv.

\subsection{The most general case}
\label{sec:gen}

In the case where both repeated scalar blocks and full blocks can be either complex or real we have to add a dynamics for real full blocks. 
This can be done following an approach analogous to the one discussed in~\cite{GL13} and exploiting a rank-$2$ property of real blocks in the matrices. In order to derive a gradient system for this general case it is necessary to add systems of differential equations for $m\times m$ real blocks $\Delta$, which
are conveniently expressed in the form  $\Delta = U Q V^{\rm T}$ with $U,V \in \R^{m \times 2}$ having orthonormal columns and a 
$2 \times 2$ orthogonal matrix $Q$. 

A full discussion of this case is omitted for conciseness. 
However, our implemented algorithm includes this case.

\section{Fast approximation of $\mu_\BP(M)$}
\label{sec:fast}

In this section, we discuss the outer algorithm for computing a lower bound of $\mu_\BP(M)$. Since the principles are the same, we treat the case of
purely complex perturbations in detail and provide a briefer discussion on the extension to the case of mixed complex/real perturbations.

\subsection{Purely complex perturbations}
\label{sec:outc}

In the following, we let $\lambda(\eps)$ denote a continuous branch of (local) maximizers for
$$\max\limits_{\l \in \Lambda_\eps^{\BPC}(M)} |\l|,$$
computed by determining the stationary points $\Delta(\eps)$
of the system of ODEs~\eqref{Fode} (or, equivalently,~\eqref{Fode2}). 
The computation of the SSV is equivalent to the smallest solution $\eps$ of the equation $|\lambda(\eps)| = 1$. 
In order to approximate this solution, we aim at computing $\eps^\star$ such that the boundary of the $\eps^\star$-spectral value set
is locally contained in the unit disk and its boundary $\partial \Lambda_{\eps^\star}^{\BPC}(M)$ is tangential to the unit circle. 
This provides a lower bound $1/\eps^\star$ for $\mu_\BPC(M)$

\begin{figure}[ht]
\centerline{
\includegraphics[scale=0.43,trim= 0mm 0.01mm 0mm 0mm]{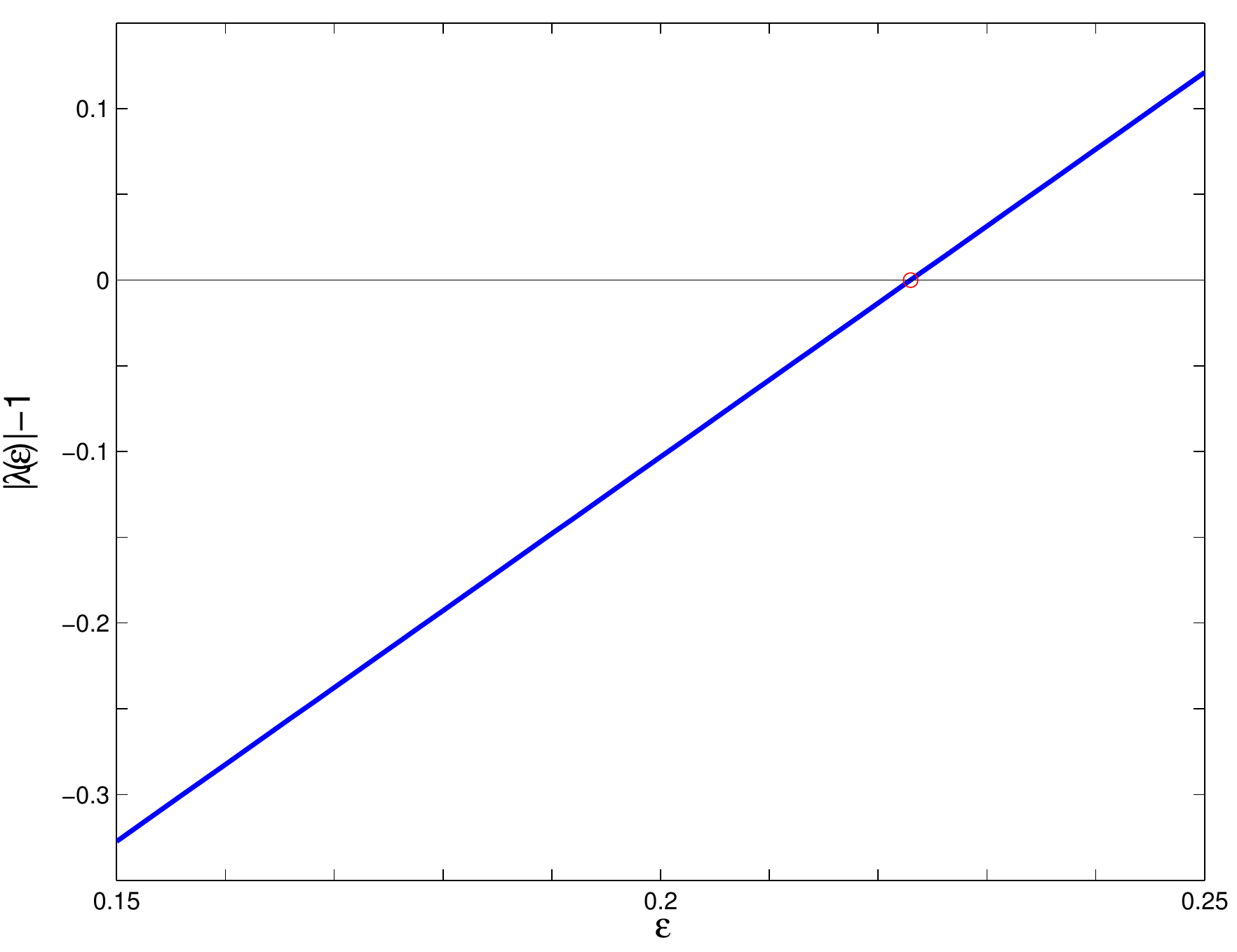}
}
\caption{The function $|\l(\eps)|-1$ for the example from Section~\ref{sec:illc}. \label{fig:lambda}}
\end{figure}

In order to apply the Newton method for solving $|\lambda(\eps)| = 1$ (see Figure~\ref{fig:lambda} for an illustration of the function $|\lambda(\eps)| - 1$), we need to compute the derivative of $|\l(\eps)|$ with respect to $\eps$.
For this purpose, we make the following generic assumption.
\begin{assumption}
For a local extremizer $\Delta(\eps)$ of $\Lambda_\eps^{\BPC}(M)$, with corresponding largest eigenvalue 
$\l(\eps)$, we assume
that $\l(\eps)$ is simple and that $\Delta(\cdot)$ and $\l(\cdot)$ are smooth in a neighborhood of $\eps$. 
\label{ass:smooth}
\end{assumption}

%


The following theorem gives an explicit and easily 
computable expression for the derivative of $|\l(\eps)|$.
\begin{theorem} \label{th:der}
Suppose that Assumption~\ref{ass:smooth} holds for $\Delta(\eps) \in \BPCO$ and $\l(\eps)$. Let
$x(\eps)$ and $y(\eps)$ be the corresponding right and left eigenvectors of 
$\eps M \Delta(\eps)$, 
scaled according to \eqref{eq:scalyx}. Consider the partitioning~\eqref{eq:partitionxz} of $x(\eps)$, $z(\eps)=M^* y(\eps)$, and suppose that Assumptions~\eqref{eq:A1} and~\eqref{eq:A2} hold.
Then
\begin{equation}
\frac{d |\l(\eps)|}{d \eps} = {} 
\frac{1}{|y(\eps)^* x(\eps)|} \Big(
\sum\limits_{i=1}^{S} |z_i(\eps)^* x_i(\eps) | + \sum\limits_{j=1}^{F} \| z_{S+j}(\eps) \|\,\| y_{S+j}(\eps) \| \Big)
> 0.
\label{eq:der}
\end{equation}
\end{theorem}
\begin{proof}
First, we observe that
\begin{equation}
\frac{d}{d\eps} |\l(\eps)| 
 =  \frac{1}{2|\l(\eps)|}\frac{d}{d\eps}|\l(\eps)|^2 = 
\frac{1}{|\l(\eps)|} \Re\left(\overline{\l}(\eps) \l^\prime(\eps)\right),
\label{eq:dl2eps}
\end{equation}
where we let ${\phantom{a}}^\prime \equiv d/d\eps$. Plugging
$$
\l^\prime(\eps)= \frac{y(\eps)^*(M \Delta(\eps)+\eps M \Delta'(\eps)) x(\eps)}
{y(\eps)^* x(\eps)}
$$
into~\eqref{eq:dl2eps} yields
\begin{eqnarray}
\frac{d}{d\eps}|\l(\eps) | 
& = & 
\frac{1}{|\l(\eps)|} \Re \left(\overline{\l}(\eps) \; 
\frac{y(\eps)^*(M \Delta(\eps)+\eps M \Delta'(\eps)) x(\eps)}{y(\eps)^*x(\eps)}\right)
\nonumber \\ 
& = & 
\frac{1}{|\l(\eps)|} \Re \left(\frac{y(\eps)^* (M \Delta(\eps)+\eps M \Delta'(\eps)) x(\eps)} 
{|y(\eps)^*x(\eps)|\e^{-i\theta(\eps)}}
|\l(\eps)|\e^{-i\theta(\eps)}\right)
\nonumber \\
& = &  
\Re \left(\frac{\langle y(\eps) x(\eps)^*, M \Delta(\eps)+\eps M \Delta'(\eps) \rangle}
{|y(\eps)^*x(\eps)|} \right). \label{eq:firstresult}
\end{eqnarray}
We now aim to prove that the second term in the sum vanishes, that is,
\begin{eqnarray}
\Re \left( y(\eps)^* M \Delta'(\eps) x(\eps) \right) & = & 0.
\label{eq:xEprimey}
\end{eqnarray}
The maximality property of the modulus of the eigenvalue $\l(\eps)$ of $\eps M \Delta(\eps)$ 
yields 
$
\Re \left( y(\eps)^* M \Delta'(\eps) x(\eps) \right) \le 0.
$
Now suppose that for some $\eps_0$, this inequality would actually be a strict inequality.
Consider $\widetilde{\Delta}(\eps) \in \BPCO$ such that $\widetilde{\Delta}(\eps_0)=\Delta(\eps_0)$ and
$\widetilde{\Delta}^\prime(\eps_0)=-\Delta^\prime(\eps_0)$. Then, for all $\eps$ sufficiently close to $\eps_0$,~\eqref{eq:firstresult} implies that 
the corresponding largest eigenvalue $\widetilde{\lambda}(\eps)$ of
$\eps M \widetilde{\Delta}(\eps)$ satisfies $|\widetilde{\lambda}(\eps)| > |\lambda(\eps_0)|$.
This, however, contradicts the extremality of $\Delta(\eps)$ and hence~\eqref{eq:xEprimey} holds.
In turn,~\eqref{eq:firstresult} gives
\[
\frac{d}{d\eps}|\l(\eps) | =  \Re \left(\frac{\langle y(\eps) x(\eps)^*, M \Delta(\eps) \rangle}
{|y(\eps)^*x(\eps)|} \right) = 
\Re \left(\frac{\langle P_\BPC \left( z(\eps) x(\eps)^* \right), \Delta(\eps) \rangle}
{|y(\eps)^*x(\eps)|} \right).
\]
%
The expression~\eqref{eq:der} now follows from the relation $\Delta(\eps) = D(\eps) P_\BP \left( z(\eps) x(\eps)^* \right)$ established 
in Theorem~\ref{stat:RF}, where the positive diagonal matrix $D(\eps)$ is such that all blocks of $\Delta(\eps)$ have unit Frobenius norm. 
The positivity of~\eqref{eq:der} is a consequence of Assumptions~\eqref{eq:A1} and~\eqref{eq:A2}.
\end{proof}

Theorem \ref{th:der} allows us to easily realize the Newton method
\begin{equation}
\eps^{(k+1)} = \eps^{(k)} - \frac{|\lambda^{(k)}|-1} 
{d|\lambda^{(k)}|},
\label{eq:Newton}
\end{equation}
where $\lambda^{(k)} = \lambda( \eps^{(k)} )$ and $d|\lambda^{(k)}|$
is the derivative of $|\lambda( \eps )|$ at $\eps = \eps^{(k)}$ given by~\eqref{eq:der}. Note that Theorem~\ref{th:der} implies local quadratic convergence of~\eqref{eq:Newton} to $\eps^\star$, provided that the assumptions of the theorem hold for $\eps = \eps^\star$.
See Table~\ref{tab:ex1} below for the numerical confirmation.

\subsection{Mixed complex/real perturbations}
\label{sec:outmix}

\begin{figure}[ht]
\centerline{
\includegraphics[scale=0.43,trim= 0mm 0.01mm 0mm 0mm]{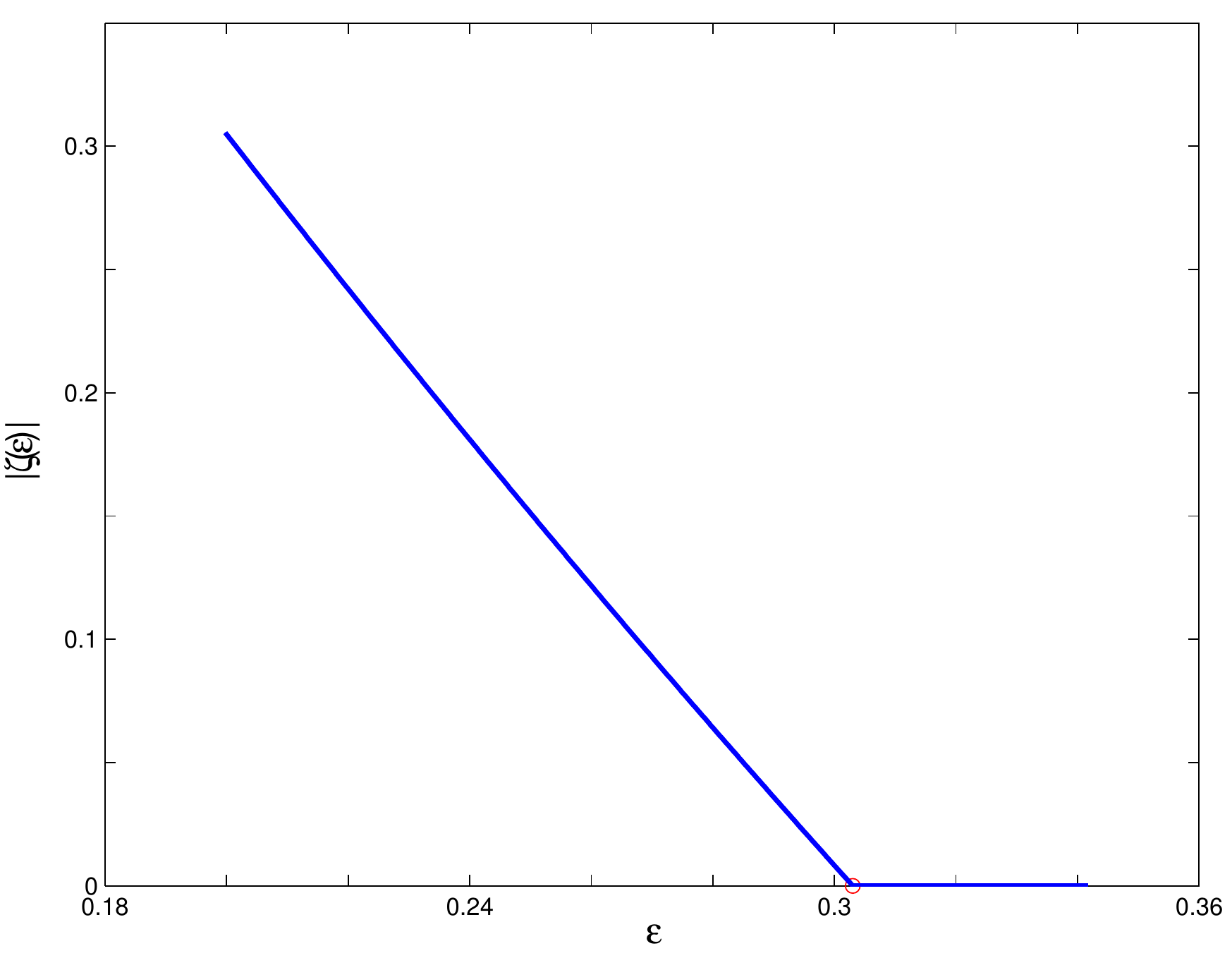}
}
\caption{The function $\z(\eps)$ for the example from Section~\ref{sec:illcr}. \label{fig:zeta}}
\end{figure}

Let $\z(\eps)$ denote a continuous branch of (local) minimizers 
of the optimization problem 
$$\min\limits_{\z \in \Sigma_\eps^{\BP}(M)} |\z|.$$ 
We aim at computing the derivative of the function $\z(\eps)$ with respect to $\eps$; see Figure~\ref{fig:zeta} for an illustration.
Since the function has a kink at the intersection with the horizontal axis,
the use of the derivative in a Newton method is meaningful only at values $\eps$ with
$\z(\eps) \neq 0$.

We make the following assumption analogous to Assumption \ref{ass:smooth}. 
\begin{assumption}
For a local extremizer $\Delta(\eps)$ of $\Sigma_\eps^{\BP}(M)$, with corresponding smallest eigenvalue 
$\zeta(\eps)$ of $\Id - \eps M \Delta(\eps)$, we assume
that $\zeta(\eps) \not=0$ is simple and that $\Delta(\cdot)$ and $\zeta(\cdot)$ are smooth in a neighborhood of $\eps$. 
\label{ass:smoothrc}
\end{assumption}

The following result is the analogue of Theorem~\ref{th:der}; it gives an explicit and easily 
computable expression for the derivative of $|\zeta(\eps)|$.
Its proof is omitted for brevity, due to its similarity with the proof of Theorem~\ref{th:der}.
\begin{theorem} \label{th:derrc}
Suppose that Assumption~\ref{ass:smoothrc} holds for $\Delta(\eps)$ and $\zeta(\eps)$. Let
$x(\eps)$ and $y(\eps)$ be the corresponding right and left eigenvectors of 
$\Id - \eps M \Delta(\eps)$, 
scaled according to \eqref{eq:scalyx}. Consider the partitioning~\eqref{eq:partitionxz} of $x(\eps)$, $z(\eps)=M^* y(\eps)$, and suppose that Assumptions~\eqref{eq:M1}--\eqref{eq:M3} hold.
Then
\begin{equation}
\begin{aligned}
\frac{d |\z(\eps)|}{d \eps} = {}- 
\frac{1}{|y(\eps)^* x(\eps)|} \biggl( &
\sum\limits_{i=1}^{S'} |z_i(\eps)^* x_i(\eps) | + \sum\limits_{i=1}^{S'} |\Re(z_i(\eps)^* x_i(\eps))| \biggr.
\\
& + \biggl. 
\sum\limits_{j=1}^{F} \| z_{S+j}(\eps) \|\,\| y_{S+j}(\eps) \| \biggr)
< 0.
\end{aligned}
\label{eq:dermix}
\end{equation}
\end{theorem}
We will make use of the following Newton method: For $\z^{(k)} = \z(\eps^{(k)}) \neq 0$,
\begin{equation}
\eps^{(k+1)} = \eps^{(k)} + \frac{|\z^{(k)}|} 
{d|\z^{(k)}|}
\label{eq:Newtonmix}
\end{equation}
where $d|\z^{(k)}|$ denotes the derivative of $|\z(\eps)|$ at $\eps =  \eps^{(k)}$, given by~\eqref{eq:dermix}. 
Note that this formula cannot be used if $0 \in \Sigma_{\eps^{(k)}}^{\BP}(M)$. 

Finally we remark that the discussed approach also allows to fix a threshold $\tau$ and approximate the problem
$$\min\limits_{\eps > 0} \mathcal{B}_\tau \cap \Sigma_{\eps}^{\BP}(M) = \emptyset,$$
where $\mathcal{B}_\tau$ is the sphere of radius $\tau$ in the complex plane. 

\subsection{Summary}
\label{sec:algoR}
Algorithm~\ref{algo-ode} describes the overall procedure for approximating the SSV
of a matrix $M$ with a prescribed block structure $\BPR$ for admissible perturbations of complex/real 
form.

\begin{algorithm}[h]
\DontPrintSemicolon
\KwData{$M$, $\BPR$, ${\rm tol} > 0$ and $\eps^{(0)},\eps_\ell$ (given lower bound), $\eps_u$ (given upper bound),
$i_{\max}$ (number of starting eigenvalues)} 
\KwResult{$\eps_f$ (approximation of $\eps^\star$)}
\Begin{
\nl \For{$i \leftarrow 1$ \KwTo $i_{\max}$}{
    Solve system of ODEs~\eqref{Fodecr2} starting from the initial matrix $\Delta_i(0)$ given by \eqref{eq:Delta0} 
		associated to the $i$-th largest eigenvalue of $M$\;
\nl Let $\Delta_i$ be the computed stationary solution and $\lz_i$
the smallest eigenvalue of $\Id - \eps^{(0)} M \Delta_i$\; 
}
\nl Set $i_{*} = \arg\min_{1 \le i \le i_{\max}} |\lz_i|$\;
\nl Set $\Delta^{(0)} = \Delta_{i_{*}}$,
$\lz^{(0)} = \lz_{i_{*}}$, $x^{(0)}, y^{(0)}$ the associated eigenvectors\;
\nl Compute $\eps^{(1)}$ by one step of the Newton method~\eqref{eq:Newtonmix}\;
\nl Set $k=1$\;
\nl \While{$|\eps^{(k)} - \eps^{(k-1)}| \ge {\rm tol}$}{
\nl Solve ODEs \eqref{Fodecr2} with $\eps=\eps^{(k)}$, starting from $\Delta(0)=\Delta^{(k-1)}$\;
\nl Let $\Delta^{(k)}$ be the stationary solution of \eqref{Fodecr2}\;	
\nl Let $\lz^{(k)}$ be the smallest eigenvalue of $\Id - \eps^{(0)} M \Delta^{(k)}$\;	
\nl \eIf{$|\lz^{(k)}| > {\rm tol}$}{ 
    Set $\eps_\ell = \eps^{(k)}$\;
    \nl Compute $\eps^{(k+1)}$ by one step of the Newton method~\eqref{eq:Newtonmix}}
		{Set $\eps_u = \eps^{(k)}$\;
		Set $\eps^{(k+1)} = (\eps_\ell+\eps_u)/2$}
\nl Set $k=k+1$.\;		   		
		}
\nl Set $\eps_f = \eps^{(k)}$.\;
}
\caption{Basic algorithm for computing $\eps^\star$ \label{algo-ode}}
\end{algorithm}

For the numerical integration of the ODEs we have made use of the forward Euler method with the step size
controlled by the monotonicity of the extremal eigenvalue. The stopping rule is based on 
two criteria, the first is the condition for the stepsize to not decrease under a prescribed minimal
value and the second relies on the difference of the extremal eigenvalues in two subsequent steps, which
should not decrease under a given tolerance. More sophisticated numerical integrators might be the
object of future research.
As for the value of initially tested eigenvalues, in our implementation of Algorithm \ref{algo-ode} 
we have made the choice  $i_{\max} = \max(n/5,5)$ for problems of dimension $n \ge 5$ and $i_{\max}=n$
otherwise.

\section{Computational results}
\label{sec:compres}

In this section we first provide some numerical tests for small matrices and then some statistics 
on the comparison between Algorithm~\ref{algo-ode} and the classical algorithm 
implemented in the \matlab Control Toolbox \mussv{} on a larger number of matrices having size 
between $5$ to $100$. For this purpose, we have developed a prototype \matlab implementation
of Algorithm~\ref{algo-ode}. As this implementation is not particularly optimized, we do not provide timings but focus
on the quality of the lower bounds. The fine tuning and efficient implementation of Algorithm~\ref{algo-ode} is beyond the scope
of this paper and subject to future work.

\subsection{Numerical tests}
\label{sec:numill}

In the following examples we consider real / complex perturbations of the form~\eqref{eq:BPcr}
and do not impose a particular order of appearance of repeated scalar blocks and full blocks, which has been done
for notational convenience only.

\subsection*{Example 1} Consider the following matrix from \cite{P16},
\begin{eqnarray*}
M & = & \left( \begin{array}{rrrrr}
 \iu       &       \frac12 - \frac12 \,\iu    &       1    &     1   &  \frac12 \\[1mm]                     
 \frac12       &      -\frac12              &      \iu   &    \iu  &  \frac12 -\frac12 \,\iu \\[1mm]      
 \iu       &       1 -   \frac12 \,\iu    &       1    &    \frac12  &   0  \\                 
-\frac12       &       \frac12 + \,\iu        &     -\frac12 + \frac12 \,\iu   &  1 + \frac12 \,\iu &  \frac12 - \frac12 \,\iu \\[1mm]   
 \frac12 + \,\iu &       \frac12 + \frac12 \,\iu    &       0    &   -\frac12 - \frac12 \,\iu & \frac12 - \frac12 \,\iu  
	\end{array} \right),
\end{eqnarray*}
along with the perturbation set
\[
\BPR = \big\{ \diag \left( \delta_1 \Id_{3}, \Delta_1 \right), \delta_1 \in \R, \ \Delta_1 \in \C^{2,2} \big\}.
\] 
Applying {\sc Matlab}'s {\tt mussv}, we obtain the perturbation $\widehat\eps \widehat\Delta$
with  
\begin{eqnarray*}
\widehat\Delta =\scriptsize 
\left( \begin{array}{rrrrr}
  1 & 0 & 0 & 0 & 0 \\                     
  0 & 1 & 0 & 0 & 0 \\                      
  0 & 0 & 1 & 0 & 0 \\
  0 & 0 & 0 & -0.661871043 - 0.048777846 \iu & 0.114656146 - 0.401316361 \iu \\
  0 & 0 & 0 & -0.325067916 + 0.037935543 \iu & 0.018201081 - 0.205013390 \iu
\end{array} \right).
\end{eqnarray*}
and $\widehat\eps = 0.546726635$, which gives the lower bound $\mu_\BPR(M) \ge \lbo = 1.829067647$. 
We note immediately that the $2 \times 2$ full complex block of $\widehat\Delta$ has norm $0.87\ldots$, which violates the norm-$1$ condition for all blocks of an extremizer given by Theorem \ref{th:max}. Consequently we expect to be able to improve the lower bound.
\begin{table}[h]
\caption{Values of $\eps^{(k)}$ and $|\z^{(k)}|$ computed by Algorithm~\ref{algo-ode} applied to Example 1.\label{tab:ex1}}
\begin{center}
\begin{tabular}{l|l|l|}\hline
$k$ & $\eps^{(k)}$ & $|\z^{(k)}|$  \\[0.1cm]
\hline 
\rule{0pt}{9pt}
\!\!\!\! $0$ & $0.321154624817$ & $0.325206140643$  \\
 $1$         & $0.475935094375$ & $6.509334991219 \cdot 10^{-6}$  \\
 $2$         & $0.475938192593$ & $2.803806976489 \cdot 10^{-12}$  \\
 $3$         & $0.475938192594$ & $1.037255102712 \cdot 10^{-16}$  \\[0.1cm]
 \hline
\end{tabular}
\end{center}
\end{table}

Table~\ref{tab:ex1} shows the result of Algorithm~\ref{algo-ode} with
$\eps^{(0)} = 1/\| M \|_2$.
This gives
\[
\eps^\star \approx \eps_3 \ \Longrightarrow \ \lbn \approx 2.101113160408110,
\]
which is very close to the upper bound $\mu_\BPR(M) \le 2.110047520373674$ computed by \mussv, 
and hence provides a sharp estimate.
The extremal perturbation $\eps^\star \Delta^\star$ is given by 
\begin{eqnarray*}
\Delta^\star = \scriptsize 
\left( \begin{array}{rrrrr}
  1 & 0 & 0 & 0 & 0 \\                     
  0 & 1 & 0 & 0 & 0 \\                      
  0 & 0 & 1 & 0 & 0 \\
  0 & 0 & 0 & 0.8414902738 - 0.0310321080 \iu & -0.0774400898 + 0.4310267415 \iu \\
  0 & 0 & 0 & 0.2196113059 + 0.1726644616 \iu & -0.1120620799 + 0.0924665133 \iu
\end{array} \right).
\end{eqnarray*} 
The obtained result compares favorably with the approximate value $2.1007$ computed in \cite{P16}.

\subsection*{Example 2} Consider 
\begin{eqnarray*}
M & = & \scriptsize \left( \begin{array}{rrrrrrrrrr} 
 -0.43 &  0.90 & -0.61 &  1.03 &  0.98 &  2.00 &  0.05 &  0.14 &  0.86 &  0.02 \\
 -0.17 & -1.84 & -1.22 & -0.35 & -0.30 &  0.95 &  1.75 & -1.64 &  0.11 & -0.05 \\
 -0.22 &  0.07 &  0.32 &  1.01 &  1.14 & -0.43 &  0.16 & -0.76 &  0.40 &  1.70 \\
  0.54 &  0.04 & -1.34 &  0.63 & -0.53 &  0.65 & -1.24 & -0.82 &  0.88 & -0.51 \\
  0.39 &  2.23 & -1.03 & -0.21 &  0.97 & -0.36 & -2.19 &  0.52 &  0.18 &  0.00 \\
  0.75 & -0.07 &  1.33 & -0.87 & -0.52 &  0.71 & -0.33 & -0.01 &  0.55 &  0.92 \\
  1.78 & -0.51 & -0.42 & -1.04 &  0.18 &  1.42 &  0.71 & -1.16 &  0.68 &  0.15 \\
  1.22 &  0.24 & -0.14 & -0.27 &  0.97 & -1.60 &  0.32 & -0.01 &  1.17 &  1.40 \\
 -1.28 &  0.25 &  0.90 & -0.44 & -0.41 &  1.03 &  0.41 & -0.69 &  0.48 &  1.03 \\
 -2.33 &  0.07 & -0.30 & -0.41 & -0.44 &  1.46 & -0.58 & -0.67 &  1.41 &  0.29
\end{array} \right)
\end{eqnarray*}
and 
\[
\BPR = \big\{ \diag \left( \delta_1 \Id_{1}, \delta_2 \Id_{1}, \delta_3 \Id_{1}, \delta_4 \Id_{2}, 
\Delta_1 \right), \delta_1, \delta_2 \in \R, \ \delta_3, \delta_4 \in \C, \ \Delta_1 \in \C^{5,5} \big\}.
\] 
Applying {\sc Matlab}'s {\tt mussv} gives the perturbation $\widehat\eps \widehat\Delta$
with  
\begin{eqnarray*}
&& \widehat\Delta = \scriptsize
\left( \begin{array}{rrrrrr}
 -1 &  0 &  0 & 0 & 0 & \mathbf{0}^{\rm T} \\                     
  0 & -1 &  0 & 0 & 0 & \mathbf{0}^{\rm T} \\                      
  0 &  0 & -1 & 0 & 0 & \mathbf{0}^{\rm T} \\
  0 &  0 &  0 & 1 & 0 & \mathbf{0}^{\rm T} \\
  0 &  0 &  0 & 0 & 1 & \mathbf{0}^{\rm T} \\
	\mathbf{0} & \mathbf{0} & \mathbf{0} & \mathbf{0} & \mathbf{0} & \hat{u} \hat{v}^{\rm T}  
\end{array} \right), \quad
\hat{u} = \scriptsize \left( \begin{array}{r} 
   0.93916167 \\
   0.06094908 \\
  -0.22409849 \\
   0.25285464 \\
  -0.01024501
\end{array} \right),	
\quad 
\hat{v} = \left( \begin{array}{r}
   0.21233474 \\
   0.27182946 \\
  -0.57210258 \\
   0.41515717 \\
   0.61754828  
\end{array} \right)	
\end{eqnarray*}
and $\widehat\eps = 0.23674574$, yielding the lower bound 
$\mu_\BPR(M) \ge \lbo = 4.22394088$.
The smallest eigenvalue of the matrix 
$\widehat\eps M \widehat\Delta$ is correctly $\lambda_1=1$.

 \begin{table}[h]
 \caption{Values of $\eps^{(k)}$ and $|\z^{(k)}|$ computed by Algorithm~\ref{algo-ode} applied to Example 2.\label{tab:ex2}}
\begin{center}
\begin{tabular}{l|l|l|}\hline
$k$ & $\eps^{(k)}$ & $|\z^{(k)}|$  \\[0.1cm]
\hline 
\rule{0pt}{9pt}
\!\!\!\! $0$ & $0.201123467713$ & $0.117799670760$  \\
 $1$         & $0.227979361395$ & $1.519798379258 \cdot 10^{-10}$  \\
 $2$         & $0.227979361429$ & $3.812788529246 \cdot 10^{-16}$  \\[0.1cm]
 \hline
\end{tabular}
\end{center}
\end{table} 
Algorithm~\ref{algo-ode} applied to this example results in Table \ref{tab:ex2}, with $\eps^{(0)} = 1/\| M \|_2$. 
The final (locally) extremal perturbation $\eps^\star \Delta^\star$ is given by 
\begin{eqnarray*}
\Delta^\star = \scriptsize
\left( \begin{array}{rrrrrr}
 -1 &  0 &  0 & 0 & 0 & \mathbf{0}^{\rm T} \\                     
  0 &  1 &  0 & 0 & 0 & \mathbf{0}^{\rm T} \\                      
  0 &  0 & -1 & 0 & 0 & \mathbf{0}^{\rm T} \\
  0 &  0 &  0 & -1 & 0 & \mathbf{0}^{\rm T} \\
  0 &  0 &  0 & 0 & -1 & \mathbf{0}^{\rm T} \\
	\mathbf{0} & \mathbf{0} & \mathbf{0} & \mathbf{0} & \mathbf{0} & u v^{\rm T}  
\end{array} \right), \quad
u = \left( \begin{array}{r} 
   0.85457765 \\
  -0.04668806 \\
  -0.28462457 \\
   0.41144779 \\
   0.13121292 
\end{array} \right), \quad
v = \left( \begin{array}{r} 
   0.15895464 \\
   0.22255005 \\
  -0.28570067 \\
   0.49433879 \\
   0.77408603 
\end{array} \right)
\end{eqnarray*}
and $\eps^\star \approx \eps_2$. The corresponding lower bound for the $\mu$-value is
$\mu_\BPR(M) \ge \lbn =  4.38636196596$, which improves the bound $\lbo$ by about $3\%$.
Note that the upper bound computed by \mussv{} is $\mu_\BPR(M) \le 4.45340809652$.

\subsection{Numerical statistics}
\label{sec:testat}


We now consider a test set of $100$ matrices with random entries and perturbations with
randomly chosen prescribed structure. Table~\ref{tab:1} and Table~\ref{tab:2} report the obtained 
results. 

\begin{table}[h]
\caption{Comparison between Algorithm~\ref{algo-ode} and {\sc Matlab}'s \mussv{}. The size of the randomly generated examples is given in the first column.
The second column shows the number of cases (among a total number of $100$) where the lower bound $\lbn$ computed with Algorithm~\ref{algo-ode} and the lower bound $\lbo$ computed 
by \mussv{} are equal, within a tolerance $10^{-3}$. Third column shows
the number of cases where Algorithm~\ref{algo-ode} is better than \mussv{} and in the fourth 
column the number of cases  where the opposite holds.} 
\begin{center}
\begin{tabular}{||*{10}{c||}||} 
\hline
\bf $n$ & \bfseries $\lbn = \lbo$ &\bfseries $\lbn > \lbo$  & \bfseries $\lbn < \lbo$ \\
\hline
$5$    & $63$ & $26$ & $11$  \\	
$10$	 & $66$ & $24$ & $10$  \\
$25$ 	 & $39$ & $50$ & $11$  \\ 
$50$ 	 & $37$ & $57$ &  $6$  \\
$100$  & $34$ & $63$ &  $3$  \\
\hline
\end{tabular}
\end{center}
\label{tab:1} 
\end{table}

\begin{table}[h]
\caption{Statistics on the difference $\delta = \lbn-\lbo$, between the lower bound $\lbn$ computed by Algorithm~\ref{algo-ode}
and the lower bound $\lbo$ computed by \mussv{}.
The second column shows the maximal difference (i.e., in favor of $\lbn$) and the third 
column shows the minimal difference (i.e., in favor of $\lbo$). The fourth and fifth column show the 
computed mean and variance, respectively.} 
\begin{center}
\begin{tabular}{||*{10}{c||}||} 
\hline
\bf $n$ & \bfseries $\delta_{\max}$ & \bfseries $\delta_{\min}$ 
        & \bfseries $\langle \delta \rangle$ &\bfseries ${\rm var}(\delta)$ \\
\hline
$5$    & $0.8115$  & $-1.1650$ & $0.0277$ & $0.0454$ \\ 
$10$	 & $1.0082$  & $-0.8674$ & $0.0364$ & $0.0382$ \\ 
$25$ 	 & $1.6358$  & $-0.5046$ & $0.1506$ & $0.1070$ \\ 
$50$ 	 & $0.8504$  & $-0.0016$ & $0.1775$ & $0.0586$ \\ 
$100$  & $6.1290$  & $-0.0782$ & $0.5793$ & $1.2956$ \\ 
\hline
\end{tabular}
\end{center}
\label{tab:2} 
\end{table}


For sizes $n=25$, $n=50$ and $n=100$, our new method performs significantly
better (that is, beyond the tolerance $10^{-3}$) in more than the half of the cases compared to the \matlab Control
Toolbox. 

Finally, let us mention the trivial consideration that one can always take the maximum of the lower bounds
by Algorithm~\ref{algo-ode} and \mussv{}. A little less trivial, one can take the output of \mussv{} to initialize Algorithm~\ref{algo-ode}, see Section~\ref{sec:delta0}.
Any such hybrid algorithm will improve upon \mussv{} and, as Tables~\ref{tab:1} and~\ref{tab:2} show, this improvement can often be quite significant. We therefore propose to complement \mussv{} with such a hybrid strategy.

\subsection{A possible combination with \mussv}

A possible combination of \mussv{} with the method presented in this article could be as follows. Whenever there appear blocks of norm smaller than one in the normalized extremizer
computed by \mussv{}, one can reduce the value $\widehat\eps$ and apply the numerical integrator
to the system of ODEs. The initial perturbation is chosen as the one computed by \mussv{}. One diminishes $\eps$ until the smallest eigenvalue 
of the matrix $\Id - \eps M \Delta$ is non zero. After following such a path (in $\eps$) it would be natural to make use of a few steps of Algorithm~\ref{algo-ode}.

Consider the following illustrative example:
\begin{eqnarray*}
M \! = \! \scriptsize \left( \begin{array}{rrrrrrrrrr} 
  -1 + \iu    &    0       & -1 - 2\,\iu &  -1          &   1          &      - 2\,\iu &  1 + \iu    &  1         &    0       &  2 - \iu \\
       \iu    &    1       &  1          &   1          &  -1 + \iu    &    1    + \iu & -1 + \iu    &  1         &      - \iu &  2       \\   
       \iu    &    0       &  0          &       \iu    &   0          &      - 2\,\iu & -1 + \iu    & -1 + \iu   &   -1 - \iu & -2 - \iu \\
       \iu    &   -4       &  0          &   1          &   1          &   -1 +    \iu & -1 + 2\,\iu &    - \iu   &     2\,\iu &  3 - \iu \\
   0          &      - \iu & -1 + \iu    &    2\,\iu    &  -1 + 2\,\iu &   -2 + 2\,\iu &  1 + \iu    &  2 - \iu   &    1 + \iu &  1 - \iu \\
  -2          &      - \iu &  1 + \iu    &  -1 - \iu    &     - \iu    &      - 2\,\iu &    - \iu    & -1 - \iu   &   -1 - \iu &  0       \\    
   1          &    1 - \iu &  1 - \iu    &   0          &  -1 - \iu    &   -1          &  0          &  1         &      - \iu &   2\,\iu \\
  -1          &     2\,\iu & -2 + \iu    &   1          &   1 - \iu    &    1          &  0          &  1 + \iu   &   - 2\,\iu &  1 - \iu \\
  -1 - 2\,\iu &      - \iu & -1 + \iu    &  -1 - 2\,\iu &       \iu    &    0          & -1 - \iu    &  0         &    1       &      \iu \\
     - 2\,\iu &    0       &  1 + \iu    &  -1 + \iu    &     - \iu    &    0          &      \iu    & -2 - \iu   &    0       &      \iu
\end{array} \right)
\end{eqnarray*}
with
\[
\BPR = \big\{ \diag \left( \Delta_1 , \delta_1 \Id_{4}, \delta_2 \Id_{4} \right), 
\Delta_1 \in \C^{2,2}, \delta_1, \delta_2 \in \R  \big\}.
\] 
Applying \mussv{} gives the following estimate:
\[
1.87690862\ldots \le \mu_\BPR(M) \le 5.26766965\ldots
\]
that is a significant gap.
The perturbation associated to the lower bound is $\widehat\eps \widehat\Delta$
with 
\begin{eqnarray*}
\widehat\Delta =\scriptsize 
\left( \begin{array}{cccc}
\ 0.01622800 - 0.44875053 \iu  &  \ 0.33074886 - 0.68259094 \iu & \mathbf{0}^{\rm T}  & \mathbf{0}^{\rm T} \\
-0.10388720 - 0.21700229 \iu  & -0.01277064 - 0.40618809 \iu    & \mathbf{0}^{\rm T}  & \mathbf{0}^{\rm T} \\
   \mathbf{0}                 & \mathbf{0}                      & \widehat{\delta}_1 \Id_4 & \mathbf{O}  \\                      
   \mathbf{0}                 & \mathbf{0}                      & \mathbf{O}               & \widehat{\delta}_2 \Id_4                      
\end{array} \right),
\end{eqnarray*}
where $\widehat\Delta$ has unit norm, but $\widehat{\delta}_1 = 0.37144260\ldots$ and $\widehat{\delta}_2 = -0.25823740\ldots$.
This suggests that the necessary optimality conditions are not fulfilled. Indeed this can be checked by computing the left and
right eigenvectors to the eigenvalue $1$ of $\eps M \Delta$.

Starting from the value $\eps = 0.532790989$ which is the reciprocal of the lower bound computed by \mussv, we
proceed a few steps backward and diminish $\eps$ until reaching the value $\eps_0 = 0.23$, for which we compute the point
$$
z =  -0.02016427 - 0.00149021 \iu, \qquad |z| = 0.0202192609,
$$
which locally minimizes the modulus of $\Sigma_{\eps_0}^{\BP}(M)$.

Performing three iterations of Algorithm~\ref{algo-ode} determines the final (locally) extremal perturbation 
$\eps^\star \Delta^\star$ with $\eps^* = 0.23478601$ and 
\begin{eqnarray*}
\Delta^* =\scriptsize 
\left( \begin{array}{cccc}
 0.44211256 - 0.19582232 \iu  & \ 0.38904261 - 0.75366740 \iu & \mathbf{0}^{\rm T} & \mathbf{0}^{\rm T} \\
 0.04777399 - 0.09593068 \iu  & -0.04015431 - 0.18364087 \iu  & \mathbf{0}^{\rm T} & \mathbf{0}^{\rm T} \\
   \mathbf{0}                 & \mathbf{0}                    & {}-\Id_4           & \mathbf{O}  \\                     
   \mathbf{0}                 & \mathbf{0}                    & \mathbf{O} & {}-\Id_4                      
\end{array} \right).
\end{eqnarray*}
This corresponds to the bound $4.259161456 \le \mu_\BPR(M)$, which improves significantly the one computed by \mussv{}
 with default parameters.



\section*{Conclusions}

In this article we have considered the problem of approximating structured singular values, which play an important role in robust control. Our main
results provide a characterization of extremizers and gradient systems, which can be integrated
numerically in order to provide approximations from below to the structured singular value of a matrix
subject to general complex/real block perturbations. The experimental results show the effectiveness of
the proposed method when compared to some classical algorithms proposed in the literature and implemented
in the \matlab Robust Control Toolbox.

\section*{Acknowledgments}

The first and third authors thank Christian Lubich for inspiring discussions.
The first author thanks the Italian INdAM GNCS for financial support as well as 
the center of excellence DEWS (L'Aquila). 
This work has been initiated during a Research in Pairs stay at the Mathematisches Forschungsinstitut
Oberwolfach.

\bibliographystyle{siam}
\bibliography{block}
\end{document}